\documentclass[12pt]{article}

\usepackage{amsmath,amsthm}

\parskip=\medskipamount
\def\a{\alpha}
\def\b{\beta}

\def\onehalf{{\textstyle\frac12}}

\def\cinfty#1{C^{\scriptscriptstyle\infty}(#1)}
\def\vectorfields#1{{\cal X}(#1)}

\def\hook{\mathop{\hbox to 6pt{\hrulefill}
                      \hbox{\vrule\phantom{\vbox to 7pt{}}}}}

\def\G{\Gamma}
\def\g{\gamma}

\def\fpd#1#2{\frac{\partial #1}{\partial #2}}
\def\R{{\rm I\kern-.20em R}}
\def\sode{{\sc Sode}}

\def\proof{{\sc Proof:}}
\newtheorem{prop}{\bf Proposition}
\newtheorem{lemma}{\bf Lemma}

\newtheorem{dfn}{\bf Definition}

\def\r{\rho}
\def\l{\lambda}
\def\s{\sigma}
\def\e{\eta}
\def\z{\zeta}
\def\le{J^1_{\l}E}
\def\lv{V_{\r}E}
\def\otau{\tau^0_1}
\def\otaupi{(\tau \circ \pi)^0_1}
\def\ov#1{\overline{#1}}
\def\ovotau{\ov{\tau}^0_1}
\def\JE{J^1E}
\def\JM{J^1M}
\def\pmb{\boldsymbol}
\def\pmbsig{\pmb{\sigma}}
\def\pmbeta{\pmb{\eta}}
\def\pmbxi{\pmb{\xi}}

\def\pmbzeta#1{\pmb{\zeta}_{#1}}
\def\pmbtheta{\pmb{\theta}}
\def\pmbom{\pmb{\omega}}
\def\pmbrho{\pmb{\rho}}
\def\pmba{\pmb{\a}}
\def\pmbb{\pmb{\b}}
\def\pmbg{\pmb{\g}}
\def\base#1{\pmb{e}_{#1}}

\def\baseX#1{\pmb{\mathcal X}_{#1}}
\def\baseV#1{\pmb{\mathcal V}_{#1}}
\def\forms#1{{\textstyle\bigwedge}^{#1}(\pi^\dag)}

\newcommand{\mybox}[1]{\makebox(0,0){\footnotesize{#1}}}
\newlength{\savelen}

\begin{document}

\title{Lie algebroid structures on a class of affine bundles}

\author{W.\ Sarlet\footnote{email: willy.sarlet@rug.ac.be} , T.\ Mestdag\\
{\small Department of Mathematical Physics and Astronomy }\\
{\small Ghent University, Krijgslaan 281, B-9000 Ghent, Belgium}\\[2mm]
E.\ Mart\'{\i}nez\\
{\small Departamento de Matem\'atica Aplicada}\\
{\small Universidad de Zaragoza, Mar\'{\i}a de Luna 3, E-50015 Zaragoza, Spain}}
\date{}
\maketitle

\vspace{1cm}

Keywords: Lie algebroids, affine bundles, exterior calculus.

AMS classification: 37J99, 17B60, 53C15

\vspace{2cm}

\begin{quote}
{\bf Abstract.} {\small We introduce the notion of a Lie algebroid structure on an affine
bundle whose base manifold is fibred over $\R$. It is argued that this is the framework
which one needs for coming to a time-dependent generalization of the theory of
Lagrangian systems on Lie algebroids. An extensive discussion is given of a way one can
think of forms acting on sections of the affine bundle. It is further shown that the affine
Lie algebroid structure gives rise to a coboundary operator on such forms. The concept of
admissible curves and dynamical systems whose integral curves are admissible, brings an
associated affine bundle into the picture, on which one can define in a natural way a
prolongation of the original affine Lie algebroid structure.}
\end{quote}

\newpage

\section{Introduction}

There has been a lot of interest recently in the study of dynamical systems
which have a Lie algebroid as carrying space (see e.g.\
\cite{Carin,CarMar,Clem,Liber,Mart,Wein}). A Lie algebroid is a vector bundle
$\pi :V \rightarrow M$, which comes equipped with two operators. To begin
with, there is a bracket operation on $Sec(\pi)$, the set of sections of
$\pi$, which provides it with a real Lie algebra structure . Secondly, there
is a linear bundle map $\r:V\rightarrow TM$, called the anchor map, which
establishes a Lie algebra homomorphism between $Sec(\pi)$ and the real Lie
algebra of vector fields on $M$ and does this in such a way that there is a
certain compatibility also with the module structure over $\cinfty{M}$. To be
precise, we have
\[ [\r(\s), \r(\e)] = \r ([\s, \e]) \qquad \mbox{and} \qquad
[\s , f \, \e] = f [\s, \e] + \r (\s)(f) \, \e,
\]
for all $\s, \e \in Sec(\pi)$ and $f \in \cinfty{M}$.

Weinstein's paper on Lagrangian mechanics and groupoids
\cite{Wein} roused new interest into the field of algebroids and
groupoids. Weinstein introduces `Lagrangian systems' on a Lie
algebroid by means of a Legendre-type map from $V$ to $V^*$,
associated to a given function $L$ on $V$. The local coordinate
expression of such equations reads:
\begin{equation}  \begin{array}{rcl}
\dot{x}^i &\!\!\!=\!\!\!& {\r}^i_{\a}(x) y^{\a}, \\[1mm] \displaystyle
\frac{d}{dt}\left(\fpd{L}{y^{\a}}\right) &\!\!\!=\!\!\!&
\displaystyle {\r}^i_{\a} \fpd{L}{x^i} - C_{\a \b}^{\g} y^{\b} \fpd{L}{y^{\g}},
\end{array} \label{autlagreq}
\end{equation}
where the $x^i$ are coordinates on $M$, $y^\a$ are fibre coordinates on $V$ and the
$C^\g_{\a\b}$ are structure functions coming from the Lie algebroid structure.
Applications for such model equations can be found e.g.\ in the theory of systems
with symmetries on principal fibre bundles and in rigid body dynamics.
Note that, more generally, equations of the form
\begin{equation}
\begin{array}{rcl} \dot{x}^i &\!\!\!=\!\!\!& {\r}^i_{\a}(x) y^{\a}, \\[1mm]
\dot{y}^{\a} &\!\!\!=\!\!\!& f^{\a}(x,y), \end{array}
\label{autdynsyst}
\end{equation}
were called ``second-order equations on a Lie algebroid'' by
Weinstein. They are indeed, to some extent, the analogues of
second-order dynamics on a tangent bundle. It is clear, however,
that these equations truly are second-order differential equations
only when the base manifold and the fibres have the same dimension
and $\r$ is injective. We will therefore rather call them
`pseudo-second-order ordinary differential equations',
pseudo-\sode s for short. Weinstein also raised the question
whether there would be a geometrical way of defining equations of
the form (\ref{autlagreq}), much in the line of the geometrical
construction of classical Lagrange equations, which makes use of
the intrinsic structures living on a tangent bundle.

One of us has recently resolved this issue \cite{Mart} by introducing
a kind of lifted Lie algebroid, where suitable analogues can be introduced
of the dilation vector field and the vertical endomorphism on a tangent bundle.

In the present paper, we wish to set the stage for an appropriate
generalization of this theory to non-autonomous systems of
differential equations. We believe that, for example at the level
of pseudo-second-order equations, the right generalization is not
just a matter of allowing the functions $\r^i_\a$ and $f^\a$ to
depend on time, but rather should produce equations of the form:
\begin{equation}
\begin{array}{rcl} \dot{x}^i &\!\!\!=\!\!\!& {\r}^i_{\a}(t,x) y^{\a}
+ \l^i(t,x) \\[1mm]
\dot{y}^{\a} &\!\!\!=\!\!\!& f^{\a}(t,x,y) \end{array}
\label{dynsyst}
\end{equation}
The reason for this simply is that we wish the structure of the equations to be invariant under
time-dependent coordinate transformations. As for Lagrange-type equations, our only concern at
the moment is to have an idea of what a time-dependent generalization of (\ref{autlagreq}) should
look like. Now, there is a way of developing a
kind of formal calculus of variations approach which leads to equations of the form
(\ref{autlagreq}), and in which the first set of equations are treated as constraints.
We have shown in \cite{opava} that if such an approach is adopted when the Lagrangian is allowed
to depend on time and the constraints are as in (\ref{dynsyst}), one obtains
equations of the form
\begin{equation}
\begin{array}{rcl} \dot{x}^i &\!\!\!=\!\!\!& {\r}^i_{\a}(t,x) y^{\a} + \l^i(t,x), \\[1mm]
\displaystyle\frac{d}{dt}\left(\fpd{L}{y^{\a}}\right) &\!\!\!=\!\!\!&
\displaystyle {\r}^i_{\a} \fpd{L}{x^i} - (C_{\a \b}^{\g} y^{\b} - C^{\g}_{\a})
\fpd{L}{y^{\g}} \end{array} \label{lagreq}
\end{equation}
where the functions $\r^i_{\a}$, $\l^i$, $C_{\a \b}^{\g}$,
$C^{\g}_{\a}$ satisfy the relations
\begin{eqnarray}
\r^i_{\a}\fpd{\r^j_{\b}}{x^i} - \r^i_{\b}\fpd{\r^j_{\a}}{x^i} &=&
\r^j_{\g} C_{\a \b}^{\g}, \label{oldceq}\\
\fpd{\r^j_{\b}}{t} + \l^i\fpd{\r^j_{\b}}{x^i}  - \r^i_{\b}\fpd{\l^j}{x^i}
&=&  \r^j_{\a} C^{\a}_{\b}.
\label{newceq}
\end{eqnarray}
Thus, we want to address the question of explaining the nature of the conditions
(\ref{oldceq},\ref{newceq}), which presumably should again have something to do with
a Lie algebroid structure.

Inspired by these analytical considerations, we will introduce the
notion of a Lie algebroid structure on an affine bundle
$\pi:E\rightarrow M$, where the base manifold $M$ in addition is
assumed to be fibred over $\R$. For the present paper, we will
limit ourselves to a number of basic features of such a theory. In
particular, we shall show in Section~3 that our defining relations
for an affine Lie algebroid are fully consistent with the
expectation of being able to develop an exterior differential
calculus of sections of the extended dual of this bundle (and its
exterior products). We shall further show in Section~4 that vector
fields on $E$, whose integral curves are `admissible curves' and
which in fact model differential equations of the pseudo-\sode\
type, can be identified in a natural way with special sections of
a kind of prolongation of the original affine bundle. This then
brings us to a final test, for this paper, of the internal
coherence of the newly defined structures: we will verify in
Section~5 whether the prolongation of a Lie algebroid, as
constructed in \cite{Mart} for the vector bundle situation,
carries over to the present more general situation.

The final section lists a number of other topics of interest, which will be
the subject of forthcoming publications. One of our objectives is to arrive at
an intrinsic geometrical construction of the time-dependent Lagrangian
equations of type (\ref{lagreq}). For the time being, however, the conditions
(\ref{oldceq},\ref{newceq}) merely serve as benchmarks, to be met by our model
of an affine Lie algebroid.

The basic ingredients for our theory are an affine bundle $\pi:E \rightarrow M$, where
the base manifold $M$ is further fibred over $\R$.
In Section~2, we define the concept of a Lie algebroid
on $\pi$. In Section~3, we show that the axioms for such a Lie algebroid structure give
rise to a consistent development of an exterior calculus on sections of the extended dual
of $E$. In Section~4, we discuss a special class of curves on $E$, which
are said to be admissible by the anchor map and we look into the concept of dynamical
systems whose integral curves all belong to this special class. In Section~5,
we define the prolongation $\pi_1 : \le \rightarrow E$ of $\pi : E \rightarrow M$
and show that it inherits the Lie algebroid structure from $\pi$.

\section{Affine Lie algebroids}

Let $M$ be an $(n+1)$-dimensional smooth manifold, which is fibred over $\R$,
$\tau: M \rightarrow \R$. We denote
the first jet bundle of $\tau$ by $\otau : \JM \rightarrow M$. It is an
affine bundle modelled on the bundle of tangent vectors to $M$ which are vertical
with respect to $\tau$; this vector bundle will be denoted by $\ovotau: VM \rightarrow M$.
To fix notations further, if $\pi:E\rightarrow M$ is an affine bundle and
$\ov{\pi}:V\rightarrow M$ its associated vector bundle, sections of $\pi$ will be
denoted by ordinary Greek characters, whereas boldface Greek type will be used for
sections of $\ov{\pi}$.

\begin{dfn}
An affine Lie algebroid over $M$ is an affine bundle $\pi: E \rightarrow M$,
with the following properties:
\begin{enumerate}
\item $Sec(\ov \pi)$, the set of sections of the vector bundle
$\ov{\pi}: V \rightarrow M$ on which $E$ is modelled, is equipped with a
skew-symmetric and bilinear (over $\R$) bracket $[\cdot,\cdot]$;
\item the affine space $Sec(\pi)$ acts by derivations on the real algebra
$Sec(\ov \pi)$, that is to say, if the same bracket notation $[\z,\pmbsig]$ is used to
denote the way $\z\in Sec(\pi)$ acts on $\pmbsig\in Sec(\ov \pi)$, we have
$[\z,\pmbsig]\in Sec(\ov \pi)$ and
\begin{equation}
[\z,\pmbsig_1+\pmbsig_2] \!=\! [\z,\pmbsig_1]+[\z,\pmbsig_2], \qquad
[\z + \pmbsig,\pmbeta] = [\z,\pmbeta] + [\pmbsig,\pmbeta],
\label{linearity}
\end{equation}
\begin{equation}
[\z,[\pmbsig,\pmbeta]] = [[\z,\pmbsig],\pmbeta] + [\pmbsig,[\z,\pmbeta]];
\label{formaljacobi}
\end{equation}
\item there exists an affine bundle map $\l: E \rightarrow \JM$ (over the identity on $M$),
with corresponding vector bundle homomorphism $\r: V \rightarrow VM$,
such that the following compatibility condition holds for all $f\in\cinfty{M}$,
\begin{equation}
[\z,f\pmbsig] = f\,[\z,\pmbsig] + \l(\z)(f)\pmbsig. \label{compat1}
\end{equation}
\end{enumerate}
Both the affine map $\l:E\rightarrow \JM$ and its linear part $\r:V\rightarrow VM$
will be called anchor maps.
\label{affalgdef}
\end{dfn}

Note that we make no notational distinction between, on the one hand,
the affine and linear anchor maps, regarded as maps between total spaces of bundles,
and their action on sections of bundles on the other hand. Needless to say, for the
interpretation of the bracket in the left-hand side of (\ref{compat2}), both $\l(\z)$
and $\r(\pmbsig)$ are regarded as vector fields on $M$.

Let us derive some further properties which follow from this definition.
First of all, if we replace $\z$ in (\ref{formaljacobi}) by $\z + \pmbxi$, for
an arbitrary $\pmbxi\in Sec(\ov \pi)$, it follows that $\forall \pmbxi,\pmbsig,\pmbeta
\in Sec(\ov \pi)$:
\begin{equation}
[\pmbxi,[\pmbsig,\pmbeta]]= [[\pmbxi,\pmbsig],\pmbeta] + [\pmbsig,[\pmbxi,\pmbeta]].
\label{jacobi}
\end{equation}
This means, in view of the first hypothesis, that the bracket on $Sec(\ov \pi)$
actually provides $Sec(\ov \pi)$ with a real Lie algebra structure.
Secondly, making the same substitution for $\z$ in (\ref{compat1}),
recalling that $\l(\z+\pmbxi)=\l(\z)+\r(\pmbxi)$, it follows that
\begin{equation}
[\pmbxi,f\pmbsig]= f\,[\pmbxi,\pmbsig] + \r(\pmbxi)(f)\pmbsig. \label{compat3}
\end{equation}
This means that the linear anchor map $\r:V\rightarrow VM$ defines a Lie algebra
homomorphism from $Sec(\ov \pi)$ into the real Lie algebra of vertical vector fields
on $M$, and that we have a classical Lie algebroid
structure on the vector bundle $\ov{\pi}:V\rightarrow M$ (although its image cannot
reach the whole of $TM$). Thirdly, replacing $\pmbeta$ by $f\pmbeta$ in (\ref{formaljacobi})
and making use of (\ref{compat1}) and (\ref{compat3}), one obtains the additional compatibility
property
\begin{equation}
[\l(\z), \r(\pmbsig)] = \r([\z,\pmbsig]), \label{compat2}
\end{equation}
from which it further follows that
\begin{equation}
[\r(\pmbxi),\r(\pmbsig)] = \r\left([\pmbxi,\pmbsig]\right). \label{liealghom}
\end{equation}

Remark: for an alternative and equivalent definition of an affine Lie algebroid, we
could impose first the Lie algebra structure (\ref{jacobi}) of the bracket on
$Sec(\ov \pi)$, together with the compatibility condition (\ref{compat3})
for the anchor map $\r$, and subsequently require that the properties
(\ref{linearity}-\ref{compat1}) hold true for at least one $\z\in Sec(\pi)$ and for
an affine map $\l:E\rightarrow J^1M$ whose linear part is $\r$. It then follows that
such properties hold for all $\z$.

We can now further extend the bracket operation to $Sec(\pi)$, as follows.

\begin{dfn} {\rm (i)} For $\pmbsig\in Sec(\ov \pi)$ and $\z\in Sec(\pi)$, we put
$[\pmbsig,\z] = -[\z,\pmbsig]$.

{\rm (ii)} For every two sections $\z_1, \z_2 \in Sec(\pi)$ with $\pmbzeta{12} =
\z_2 - \z_1$: $\ [\z_1,\z_2] = [\z_1, \pmbzeta{12}]$.
\end{dfn}

Observe that the extended bracket in (ii) is a map from
$Sec(\pi)\times Sec(\pi)$ to $Sec(\ov \pi)$. As we will show
below, it has Lie algebra type properties, which could justify
talking about an ``affine Lie algebra structure", were it not that
this term is in use already in the literature, with an entirely
different meaning. The extended bracket also has Lie algebroid
type properties with respect to the anchor maps $\l$ and $\r$.

\begin{prop}
The bracket $[\cdot,\cdot]: Sec(\pi)\times Sec(\pi) \rightarrow Sec(\ov \pi)$, has
the following properties:
\begin{eqnarray}
&& [\z_1,\z_2 + \pmbsig] = [\z_1, \z_2] + [\z_1, \pmbsig], \label{linearity2} \\
&& [\z_1, \z_2] = - [\z_2, \z_1], \label{skewsym} \\
&& [[\z_1,\z_2],\z_3] + [[\z_2,\z_3],\z_1]  + [[\z_3,\z_1],\z_2] = \pmb{0}, \label{jacobi2} \\
&& \r([\z_1, \z_2]) = [\l(\z_1), \l(\z_2)]. \label{compat4}
\end{eqnarray}
\end{prop}
\proof \ \ The first property follows immediately from the definition and (\ref{linearity}).
Next, we have $[\z_2,\z_1]=[\z_2,\pmbzeta{21}]=-[\z_1+\pmbzeta{12},\pmbzeta{12}]
= -[\z_1,\z_2]$. For the Jacobi identity, using a simple summation sign to indicate
the cyclic sum over the three sections in each summand, we have
\[
\sum [[\z_1,\z_2],\z_3] = \sum [[\z_1,\z_2],\z_2+\pmbzeta{23}] =
\sum [[\z_1,\z_2],\pmbzeta{23}],
\]
in view of the linearity properties and the skew-symmetry of the bracket. Substituting
subsequently $\z_1+\pmbzeta{12}$ for $\z_2$, we obtain
\[
\sum [[\z_1,\z_2],\z_3] = \sum [[\z_1,\pmbzeta{12}],\pmbzeta{23}],
\]
which is zero in view of (\ref{formaljacobi}). Finally, the compatibility property
(\ref{compat4}) easily follows in the same way from the definition of the extended
bracket and (\ref{compat2}).  $\qed$

To understand what an affine Lie algebroid structure means in
coordinates, let us coordinatize $E$ in the usual way, as follows:
$t$ denotes the coordinate  on $\R$; $(x^i)_{1\leq i\leq n}$ are
fibre coordinates on $M$; we further choose a local section $e_0$
of $\pi$ to play the role of zero section and a local basis
$({\pmb e}_\a)_{1\leq \a\leq k}$ for $Sec(\ov \pi)$. Then, if $e$
is a point in the fibre $E_m$ over $m\in M$, it can be written in
the form: $e=e_0(m) + y^\a{\pmb e}_\a(m)$; $(t,x^i,y^\a)$ are
coordinates of $e$ ($(t,x^i)$ being the coordinates of $m$).

We have
\begin{equation}
[{\pmb e}_\a,{\pmb e}_\b]= C^\g_{\a\b}(t,x){\pmb e}_\g, \qquad
[e_0,{\pmb e}_\a] = C^\b_\a(t,x) {\pmb e}_\b, \label{strucfunc}
\end{equation}
for some {\it structure functions\/} $C^\g_{\a\b}=-C^\g_{\b\a}$ and $C^\b_\a$ on $M$.
The affine map $\l$ and its linear part $\r$ are fully determined by
\begin{equation}
\l(e_0) = \fpd{}{t} + \l^i(t,x)\fpd{}{x^i}, \qquad \r({\pmb e}_\a) = \r^i_\a(t,x)\fpd{}{x^i}.
\label{anchorfunc}
\end{equation}
The further characterization of the Lie algebroid structure now has the following
coordinate translation. The derivation property (\ref{formaljacobi}) and the resulting
Jacobi identity (\ref{jacobi}) mean that we have:
\begin{eqnarray}
&&\fpd{C_{\a \b}^{\mu}}{t}
+ \l^i \fpd{C_{\a \b}^{\mu}}{x^i} + C_{\a \b}^{\g} C^{\mu}_{\g}
= C_{\a \g}^{\mu} C^{\g}_{\b} - C_{\b \g}^{\mu} C^{\g}_{\a} +
\r^i_{\a} \fpd{C^{\mu}_{\b}}{x^i} - \r^i_{\b} \fpd{C^{\mu}_{\a}}{x^i},
\label{fjaccoord} \\
&& \sum_{\a,\b,\g} \left(\r^i_{\a} \fpd{C_{\b \g}^{\mu}}{x^i}
+ C_{\a \nu}^{\mu} C_{\b \g}^{\nu}\right) = 0, \label{jaccoord}
\end{eqnarray}
where the summation this time refers to a cyclic sum over $\a,\b,\g$ and also the
compatibility conditions (\ref{compat1}),(\ref{compat3}) have been invoked. Finally,
the properties (\ref{compat2}) and (\ref{liealghom}), for which it is sufficient to
express that $[\l(e_0),\r({\pmb e}_\a)]=\r([e_0,{\pmb e}_\a])$ and
$[\r({\pmb e}_\a),\r({\pmb e}_\b)]=\r([{\pmb e}_\a,{\pmb e}_\b])$, require that
\begin{eqnarray}
&& \fpd{\r^j_{\b}}{t} + \l^i\fpd{\r^j_{\b}}{x^i}  - \r^i_{\b}\fpd{\l^j}{x^i}
=  C^{\a}_{\b}\r^j_{\a}, \label{compat2coord} \\
&& \r^i_{\a}\fpd{\r^j_{\b}}{x^i} - \r^i_{\b}\fpd{\r^j_{\a}}{x^i} =
C_{\a \b}^{\g} \r^j_{\g}. \label{liealghomcoord}
\end{eqnarray}
These are precisely the relations (\ref{oldceq}) and (\ref{newceq})
we encountered in the Introduction, in the context of Lagrangian equations of type
(\ref{lagreq}).

It is of some interest to look at the way the various structure and anchor map functions
transform under coordinate transformations. There are two distinct levels in making a
change of coordinates on $E$, which we will describe separately. Firstly, we could
choose a different (local) zero section ${\ov e}_0$ and a different local basis
${\ov{\pmb e}}_\b$ for $Sec(\ov \pi)$: say that ${\pmb e}_\a=A^\b_\a {\ov{\pmb e}}_\b$
and $e_0={\ov e}_0 + B^\a {\ov{\pmb e}}_\a$. This amounts to making an affine change of
coordinates in the fibres of the form: ${\ov y}^\a = A^\a_\b(t,x)y^\b + B^\a(t,x)$.
Putting $[{\ov{\pmb e}}_\a,{\ov{\pmb e}}_\b]= {\ov C}^\g_{\a\b}{\ov{\pmb e}}_\g$,
$[{\ov e}_0,{\ov{\pmb e}}_\a] = {\ov C}^\b_\a {\ov{\pmb e}}_\b$, and also
$\l({\ov e}_0) = \fpd{}{t} + {\ov \l}^i\fpd{}{x^i}$, $\r({\ov{\pmb e}}_\a) = {\ov \r}^j_\a
\fpd{}{x^j}$, one can verify that the following transformation rules apply:
\[
\r^i_\a=A_\a^\b{\ov \r}^i_\b, \qquad {\ov \l}^i = \l^i - B^\a{\ov \r}^i_\a,
\]
and further
\begin{eqnarray*}
C_{\a\b}^{\g} A^{\mu}_{\g} &=& {\ov C}^{\mu}_{\g\nu}A^{\g}_{\a}A^{\nu}_{\b} +
\r^i_{\a} \fpd{A_{\b}^{\mu}}{x^i} - \r^i_{\b} \fpd{A_{\a}^{\mu}}{x^i}, \\
C_{\b}^{\g} A_{\g}^{\a} &=& {\ov C}_{\mu}^{\a} A_{\b}^{\mu} +
{\ov C}_{\g\mu}^{\a} B^{\g} A_{\b}^{\mu}
+ \fpd{A_{\b}^{\a}}{t} + \l^i \fpd{A_{\b}^{\a}}{x^i}  - \r^i_{\b} \fpd{B^{\a}}{x^i}.
\end{eqnarray*}
At a different level, one can make a change of coordinates on $M$, of the form:
$t'=t,\ {x'}^i={x'}^i(t,x)$. This has an effect on the anchor map functions of the form:
\[
{\r'}^j_\a = \r^i_\a\fpd{{x'}^j}{x^i}, \qquad {\l'}^j=\fpd{{x'}^j}{t} +
\l^i\fpd{{x'}^j}{x^i}.
\]
A general change of adapted coordinates is of course a composition of the two steps
described above.

\section{Exterior calculus on an affine Lie algebroid}

We first recall some features of the by now standard theory of Lie
algebroids on a vector bundle (see \cite{Mac}). Considering
sections of exterior powers of the dual bundle, one gets a notion
of forms on sections of the vector bundle, on which an exterior
derivative can be defined which involves the Lie algebroid bracket
and the anchor map. It then turns out that the Jacobi identity of
the Lie algebroid bracket and the compatibility with the bracket
of vector fields via the anchor map are exactly the conditions for
this exterior derivative to have the co-boundary property $d^2=0$
(see also \cite{GrabUrbI,GrabUrbII,Mart}). In our opinion, such a
feature in itself gives a strong indication that the
generalization from Lie algebra to Lie algebroid is indeed a
meaningful step. We shall therefore investigate in this
section whether a similar support can be detected for our
extension to Lie algebroids on affine bundles.

The extended dual of the affine space $E_m$ is the space of real valued affine
functions on $E_m$ and will be denoted by $E^{\dag}_m$. The union of these spaces
over all points $m\in M$ gives us a bundle $\pi^{\dag}: E^{\dag} \rightarrow M$ say.
Although this is in fact a vector bundle, we are interested in the action of its sections
(and sections of its exterior powers) on sections of the affine bundle $\pi$. This brings
some subtleties into the picture which need to be investigated in sufficient detail.
We will write ${\ov \pi}^*:V^*\rightarrow M$ for the dual bundle of ${\ov \pi}$ and also
use boldface type for its sections (and the sections of its exterior powers).
Now, to begin with, if $\theta$ is a section of $\pi^\dag$ and $\z\in Sec(\pi)$,
$\theta(\z)$ is a function on $M$ defined by $\theta(\z)(m)=\theta_m(\z_m)$. $\theta_m$
being an affine function on $E_m$, there exists an associated element $\pmbtheta_m\in V^*_m$
such that $\forall e_m\in E_m,\;\pmbsig_m\in V_m$, we have $\theta_m(e_m+\pmbsig_m)=
\theta_m(e_m)+\pmbtheta_m(\pmbsig_m)$. Expressed in slightly different terms and now at the
level of sections again, a $\theta\in Sec(\pi^\dag)$ is such that there exist a $\theta_0
\in Sec(\pi^\dag)$ and a $\pmbtheta\in Sec({\ov \pi}^*)$, such that for all $\z\in Sec(\pi)$:
\begin{equation}
\theta(\z) = \theta_0(\z_0) + \pmbtheta(\pmbzeta{}), \label{oneform}
\end{equation}
where $\z_0$ is any section and then $\z=\z_0+\pmbzeta{}$. The two composing elements
$\theta_0$ (which in fact is simply $\theta$ itself here) and $\pmbtheta$ do not depend
on the choice of $\z_0$. With sections of $\pi^\dag$ as our notion of 1-forms on $Sec(\pi)$,
there is of course no linearity with respect to multiplication by functions on $M$. We can
now come in a similar way to the following concept of $k$-forms on $Sec(\pi)$ (thereby
taking for granted that the meaning of a $k$-form on a vector bundle such as
$Sec({\ov \pi})$ is known).

\begin{dfn} A $k$-form on the affine bundle $Sec(\pi)$ ($k\geq 1$) is a map
$\omega: Sec(\pi)\times \cdots \times Sec(\pi)\rightarrow \cinfty{M}$,
for which there exists a $k$-form $\pmbom$ on the associated vector bundle $Sec({\ov \pi})$
and a map $\omega_0: Sec(\pi)\times Sec({\ov \pi}) \times\cdots\times
Sec({\ov \pi})\rightarrow \cinfty{M}$ with the following properties:
\begin{enumerate}
\item $\omega_0$ is skew-symmetric and $\cinfty{M}$-linear in its $k-1$ vector arguments;
\item $\forall \z\in Sec(\pi)$ and $\forall\pmbsig,\pmbzeta{j}\in Sec({\ov \pi})$, we have
\begin{equation}
\omega_0(\z+\pmbsig,\pmbzeta{1},\ldots,\pmbzeta{k-1}) =
\omega_0(\z,\pmbzeta{1},\ldots,\pmbzeta{k-1}) +
\pmbom(\pmbsig,\pmbzeta{1},\ldots,\pmbzeta{k-1}); \label{ompmbom}
\end{equation}
\item $\forall \z_i\in Sec(\pi)$, if we choose an arbitrary $\z_0\in Sec(\pi)$ and put
$\z_i=\z_0+\pmbzeta{i}$, we have
\begin{equation}
\omega(\z_1,\ldots,\z_k) = \sum_{i=1}^k(-1)^{i-1}\omega_0(\z_0,\pmbzeta{1},\ldots,
\hat{\pmbzeta{i}},\ldots,\pmbzeta{k})
+ \pmbom(\pmbzeta{1},\ldots,\pmbzeta{k}). \label{kform}
\end{equation}
\end{enumerate}
\end{dfn}

There are a number of properties to be checked to make sure that this definition makes
sense. First of all, one can verify that with two different choices of a reference
section, $\z_0$ and $\z'_0$ for example, related through $\z_0=\z'_0+\pmbsig$, it follows
from the second requirement that
\begin{eqnarray*}
\lefteqn{\sum_{i=1}^k(-1)^{i-1}\omega_0(\z_0,\pmbzeta{1},\ldots,
\hat{\pmbzeta{i}},\ldots,\pmbzeta{k}) + \pmbom(\pmbzeta{1},\ldots,\pmbzeta{k}) }\\
&&=\; \sum_{i=1}^k(-1)^{i-1}\omega_0(\z'_0,\pmbzeta{1}',\ldots,
\hat{\pmbzeta{i}'},\ldots,\pmbzeta{k}') + \pmbom(\pmbzeta{1}',\ldots,\pmbzeta{k}').
\end{eqnarray*}
Secondly, the two elements $\omega_0$ and $\pmbom$ which make up $\omega$ are unique.
Indeed, assuming there would be a second couple $\omega'_0$ and $\pmbom'$ making up
the same $\omega$, it follows by choosing $\z_0=\z_1$ (such that $\pmbzeta{1}=0$) that
$\omega_0=\omega'_0$, after which it is clear that also $\pmbom=\pmbom'$.
Note finally that the definition implies that $\omega$ itself is skew-symmetric in all
its arguments.

The set of forms on $Sec(\pi)$, which we will denote by $\forms{}$,
is a module over the ring $\cinfty{M}$, which also constitutes
the set of 0-forms. The wedge product of two forms is defined in the usual way.
By way of example, if $\a$ and $\b$ are 1-forms, we have
\[
(\a\wedge\b)(\z_1,\z_2) = \a(\z_1)\b(\z_2)-\a(\z_2)\b(\z_1),
\]
which for every choice of a reference section $\z_0$ gives rise to:
\begin{eqnarray}
(\a\wedge\b)(\z_1,\z_2) &=& \a(\z_0)\pmbb(\pmbzeta{2}-\pmbzeta{1}) -
\b(\z_0)\pmba(\pmbzeta{2}-\pmbzeta{1}) \nonumber \\
&& \mbox{} + (\pmba\wedge\pmbb)(\pmbzeta{1},\pmbzeta{2}). \label{2wedge}
\end{eqnarray}
It follows that $\pmba\wedge\pmbb$ is the 2-form on $Sec({\ov \pi})$ corresponding to
$\a\wedge\b$, and
\begin{equation}
(\a\wedge\b)_0(\z,\pmbsig)= \a(\z)\pmbb(\pmbsig) - \b(\z)\pmba(\pmbsig). \label{2wedge-0}
\end{equation}
Similarly, for the wedge product of three 1-forms, we have
\begin{equation}
(\a\wedge\b\wedge\g)_0(\z,\pmbzeta{1},\pmbzeta{2}) =
\Big(\a(\z)(\pmbb\wedge\pmbg)+\b(\z)(\pmbg\wedge\pmba)+\g(\z)(\pmba\wedge\pmbb)\Big)
(\pmbzeta{1},\pmbzeta{2}). \label{3wedge-0}
\end{equation}
These examples suggest to formalize the representation of $k$-forms a bit
further. As a preliminary remark, it may sometimes be of interest to extend the
interpretation of the operator $\omega_0$ in such a way that its single affine
section argument need not necessarily be the first. This
can simply be achieved by declaring $\omega_0$ to be
skew-symmetric in all its arguments (but still $\cinfty{M}$-linear
in its vector arguments only). More importantly, we shall take the sum of
$\omega_0$-terms in (\ref{kform}) to define another operator, denoted by $\omega^0$,
as follows:
\begin{equation}
\omega^0(\z_1,\ldots,\z_k) = \sum_{i=1}^k(-1)^{i-1}\omega_0(\z_0,\pmbzeta{1},\ldots,
\hat{\pmbzeta{i}},\ldots,\pmbzeta{k}) \;=\;
\sum_{i=1}^k \omega_0(\pmbzeta{1},\ldots,
\underset{i}{\check{\z_0}},\ldots,\pmbzeta{k}), \label{omega^0}
\end{equation}
where the second expression takes the above remark into account and
the symbol $\underset{i}{\check{\z_0}}$ then indicates that $\z_0$ has been inserted
in the $i$-th argument. The other important new convention we will adopt is to regard
$\pmbom$ also as acting on affine sections:
\begin{equation}
\pmbom(\z_1,\ldots,\z_k)=\pmbom(\pmbzeta{1},\ldots,\pmbzeta{k}). \label{newpmbom}
\end{equation}
This way, we can formally write
\begin{equation}
\omega=\omega^0 + \pmbom, \label{newomega}
\end{equation}
whereby it is to be understood that the two composing terms $\omega^0$ and $\pmbom$
are not $k$-forms on $Sec(\pi)$ by themselves. In fact, to compute their value when acting
on $k$ sections $\z_i$, a reference section $\z_0$ has to be chosen, but as argued above,
the value of the sum $\omega^0 + \pmbom$ in the end does not depend on that choice.

The rather formal looking decomposition (\ref{newomega}) now greatly facilitates the
representation of wedge products and will make the general coordinate representation of a
form more transparant. For example, the result (\ref{3wedge-0}) means that
\begin{equation}
(\a\wedge\b\wedge\g)_0=\a^0\otimes(\pmbb\wedge\pmbg) + \b^0\otimes(\pmbg\wedge\pmba)
+ \g^0\otimes(\pmba\wedge\pmbb), \label{awbwg_0}
\end{equation}
which then implies from (\ref{omega^0}) that
\begin{equation}
(\a\wedge\b\wedge\g)^0=\a^0\wedge\pmbb\wedge\pmbg + \pmba\wedge\b^0\wedge\pmbg
+ \pmba\wedge\pmbb\wedge\g^0, \label{awbwg^0}
\end{equation}
as expected. More generally, it follows directly from the defining formula for wedge products
that for $\omega=\omega^0 + \pmbom$ and $\rho=\rho^0+\pmbrho$:
\begin{equation}
\omega\wedge\rho = \omega^0\wedge\rho^0+\omega^0\wedge\pmbrho + \pmbom\wedge\rho^0
+ \pmbom\wedge\pmbrho, \label{wedgecalcul}
\end{equation}
where the sum of the first three terms is $(\omega\wedge\rho)^0$.

Suppose that, for a coordinatization of $E$,
we have chosen a zero section $e_0$ and a local basis of vector sections
${\pmb e}_\a$. Denote by $\{{\pmb e}^\b\}$ the dual basis for $Sec({\ov \pi}^*)$. There
exists a global section of $\pi^\dag$ which for each $m$ selects in $E^\dag_m$ the constant
function 1. We will call it $e^0$. The local
zero section $e_0$ of $Sec(\pi)$ can now play the role of the reference section $\z_0$
in our general considerations. Writing $\z=e_0+\z^\a{\pmb e}_\a$ for an arbitrary section
$\z$, we have for each 1-form $\theta$: $\theta(\z)=\theta(e_0)+\z^\a\pmbtheta({\pmb e}_\a)$.
Putting $\pmbtheta({\pmb e}_\a)=\theta_\a$ and $\theta(e_0)=\theta_0$, we see that $\theta$
has the local representation
\begin{equation}
\theta=\theta_0 e^0 + \theta_\a {\pmb e}^\a, \label{localtheta}
\end{equation}
where, in agreement with the general decomposition (\ref{newomega}), ${\pmb e}^\a$ has to be
regarded now as acting on $Sec(\pi)$ and $\theta^0=\theta_0 e^0$.
To be precise, putting
\begin{equation}
e_\a=e_0 + {\pmb e}_\a, \label{affinebasis}
\end{equation}
the action of ${\pmb e}^\b$ on affine sections, which can be given a meaning only after introducing
a reference section, is determined by:
\begin{equation}
{\pmb e}^\b(e_0)=0, \qquad {\pmb e}^\b(e_\a)=\delta^\b_\a. \label{duality}
\end{equation}
There is another slight abuse of notation in (\ref{localtheta}) since $\theta_0$ could have
a double meaning: in (\ref{localtheta}) it represents a local function on $M$, whereas it
also could refer to the operator introduced in (\ref{oneform}) and more generally in Definition~3.
We will, however, seldom use the notation $\theta_0$ in the latter sense when dealing with
coordinate calculations, so that the meaning will always be clear from the context.

Let us now see how all these notations fit together when we start wedging 1-forms.
For two 1-forms $\a=\a^0+\pmba$ and $\b=\b^0+\pmbb$ we find, for example from (\ref{2wedge-0})
and (\ref{omega^0}), that
\begin{equation}
(\a\wedge\b)^0=\a^0\wedge\pmbb - \b^0\wedge\pmba. \label{awb^0}
\end{equation}
This is in agreement with the general formula (\ref{wedgecalcul}) since obviously $\a^0\wedge\b^0=0$.
Expressing $\a$ and $\b$ with respect to the basis $(e^0,{\pmb e}^\a)$, we find
\begin{equation}
\a\wedge\b = (\a_0\b_\g-\b_0\a_\g)\,e^0\wedge{\pmb e}^\g + \onehalf(\a_\g\b_\delta
-\a_\delta\b_\g)\,{\pmb e}^\g\wedge{\pmb e}^\delta. \label{basis2form}
\end{equation}
Similarly, for the wedge product of three 1-forms with local representations of the form
(\ref{localtheta}), we obtain
\begin{eqnarray}
\a\wedge\b\wedge\g &=& \onehalf\Big(\a_0(\b_\mu\g_\nu-\b_\nu\g_\mu) +
\b_0(\g_\mu\a_\nu-\g_\nu\a_\mu) \nonumber \\
&& \mbox{} + \g_0(\a_\mu\b_\nu-\a_\nu\b_\mu)\Big)e^0\wedge{\pmb e}^\mu\wedge{\pmb e}^\nu
+ \pmba\wedge\pmbb\wedge\pmbg.
\label{basis3form}
\end{eqnarray}
It should now be clear without going into any further detail that a general $k$-form on
$Sec(\pi)$ locally has the following representation,
\begin{equation}
\omega = \frac{1}{(k-1)!}\, \omega_{0\mu_1\cdots\mu_{k-1}}\,e^0\wedge{\pmb e}^{\mu_1}\wedge
\cdots\wedge{\pmb e}^{\mu_{k-1}} +
\frac{1}{k!}\,\omega_{\mu_1\cdots\mu_{k}}\,{\pmb e}^{\mu_1}\wedge
\cdots\wedge{\pmb e}^{\mu_{k}}, \label{basiskform}
\end{equation}
where the coefficients are functions on $M$, which are skew-symmetric in all their
indices (including the zero for the first term); we have
\begin{eqnarray}
\omega_{0\mu_1\cdots\mu_{k-1}}&=& \omega(e_0,e_{\mu_1},\ldots,e_{\mu_{k-1}}), \label{compom0}\\
\omega_{\mu_1\cdots\mu_{k}} &=& \omega(e_{\mu_1},\ldots,e_{\mu_{k}})
- \sum_{i=1}^k \omega_{\mu_1\cdots\underset{i}{\check{0}}\cdots\mu_k}, \label{compom}
\end{eqnarray}
where $\underset{i}{\check{0}}$ again means that the index $\mu_i$ has been replaced by $0$.

Before arriving at our main goal, the development of an exterior calculus on forms, we will
recall a few generalities about derivations. Derivations on $\forms{}$ are defined in the usual way.
Following the standard work of Fr\"olicher and Nijenhuis \cite{FrNi}, one easily shows that
derivations are local operators and that they are completely determined by their action on
functions and 1-forms. The commutator of two derivations $D_i$, of degree $r_i$ say, is again
a derivation, of degree $r_1+r_2$, defined by
\begin{equation}
[D_1,D_2] = D_1\circ D_2 - (-1)^{r_1r_2}D_2\circ D_1. \label{commutator}
\end{equation}
Perhaps the simplest type of derivation is contraction with a section.

\begin{dfn} For $\omega\in\forms{k}$ and $\z\in Sec(\pi)$, $i_\z\omega\in\forms{k-1}$
is defined by
\begin{equation}
i_{\z} \omega (\z_1,\ldots,\z_{k-1}) = \omega(\z,\z_1, \ldots \z_{k-1}). \label{contract}
\end{equation}
\end{dfn}
The proof that this is a derivation of degree $-1$ is standard and does not depend
on the pecularities of our present theory. But perhaps we have to convince ourselves in the
first place that $i_\z\omega$ is indeed a form in the sense of Definition~3.

\begin{prop}
$i_{\z}\omega$ is a $(k-1)$-form which, in the sense of the general defining relation (\ref{kform}),
is determined by an operator $(i_{\z}\omega)_0$ and a $k$-form ${\mathbf i}_{\pmbzeta{}}\pmbom$
on $Sec(\ov \pi)$, defined as follows: for all $\pmbzeta{i}\in Sec(\ov \pi)$, $\z_0\in Sec(\pi)$,
\begin{eqnarray}
(i_{\z}\omega)_0 (\z_0,\pmbzeta 2, \ldots \pmbzeta{k-1}) &=&
- \omega_0(\z_0,\pmbzeta{},\pmbzeta 2, \ldots \pmbzeta{k-1}),
\quad \mbox{where\ \ } \pmbzeta{}=\z-\z_0, \label{contract0} \\
{\mathbf i}_{\pmbzeta{}}\pmbom(\pmbzeta 2, \ldots, \pmbzeta k) &=&
\omega_0(\z, \pmbzeta 2, \ldots, \pmbzeta k). \label{contractbold}
\end{eqnarray}
We further have the property (with $\pmbzeta{1}=\pmbzeta{}$):
\begin{equation}
i_{\z} \omega (\z_2,\ldots,\z_{k}) = \sum_{i=1}^k(-1)^{i-1}\left(i_\z\omega_0\right)
(\pmbzeta{1},\ldots,\hat{\pmbzeta i},\ldots,\pmbzeta{k}). \label{property}
\end{equation}
\end{prop}
\proof \ \ A direct computation, using (\ref{contract}) and (\ref{kform}), gives
\begin{eqnarray*}
\lefteqn{i_{\z}\omega(\z_2,\ldots,\z_{k}) =
\omega(\z,\z_2,\ldots,\z_{k})}  \\
&&=\; \omega_0(\z_0,\pmbzeta 2,\ldots, \pmbzeta{k}) + \sum_{j=1}^{k-1}(-1)^j
\omega_0(\z_0,\pmbzeta{},\pmbzeta 2, \ldots, \hat{\pmbzeta j},\ldots, \pmbzeta{k})
+ \pmbom(\pmbzeta, \pmbzeta 2, \ldots, \pmbzeta{k}) \\
&&=\; \sum_{j=1}^{k-1}(-1)^j
\omega_0(\z_0,\pmbzeta{},\pmbzeta 2, \ldots, \hat{\pmbzeta j},\ldots, \pmbzeta{k})
+ \omega_0(\z,\pmbzeta 2,\ldots, \pmbzeta{k}),
\end{eqnarray*}
from which we are led to introduce $(i_{\z}\omega)_0$ and ${\mathbf i}_{\pmbzeta{}}\pmbom$ as in
(\ref{contract0}) and (\ref{contractbold}). It is then straightforward to verify that these two
operators are linked by a property of type (\ref{ompmbom}), so the first statement follows.
Observe that, with an obvious meaning for contraction
of the operator $\omega_0$ with $\z$, we can write: ${\mathbf i}_{\pmbzeta{}}\pmbom =
i_\z\omega_0$. The somewhat peculiar feature of the additional property is that $i_\z\omega$
can be completely computed from $i_\z\omega_0$. To prove this we again start from (\ref{kform})
to write (with $\z_1=\z$)
\[
i_{\z}\omega(\z_2,\ldots,\z_{k}) = \sum_{i=1}^{k}(-1)^{i-1}
\omega_0(\z_0,\pmbzeta{1},\ldots, \hat{\pmbzeta i},\ldots, \pmbzeta{k})
+ \pmbom(\pmbzeta 1, \ldots, \pmbzeta{k}).
\]
This time, we substitute $\z_1-\pmbzeta{1}$ for $\z_0$ and observe that the second part of the sum
involving $\omega_0$, in view of (\ref{ompmbom}) then precisely cancels the last term. \qed

Before we can arrive now at the definition of an exterior derivative
operator, we need to give a meaning also to the value of a
$k$-form $\omega$, when say its first argument is taken to be a
vector section.

\begin{dfn} If $\omega$ is a $k$-form on $Sec(\pi)$, then for $\pmbsig\in Sec({\ov \pi})$
and $\z_i\in Sec(\pi)$, we put
\begin{equation}
\omega(\pmbsig,\z_2,\ldots,\z_k)=\omega(\z_1+\pmbsig,\z_2,\ldots,\z_k)
- \omega(\z_1,\z_2,\ldots,\z_k), \label{kformvector}
\end{equation}
where $\z_1$ is chosen arbitrarily.
\end{dfn}

For this to make sense of course, we need to be sure that the
result does not depend on the choice of $\z_1$. Now, if we
evaluate the right-hand side of the defining relation by using
(\ref{kform}), we obtain
\begin{equation}
\omega(\pmbsig,\z_2,\ldots,\z_k)=\pmbom(\pmbsig,\pmbzeta{2},\ldots,\pmbzeta{k})
+
\sum_{i=2}^k(-1)^{i-1}\omega_0(\z_0,\pmbsig,\pmbzeta{2},\ldots,\hat{\pmbzeta{i}},
\ldots,\pmbzeta{k}). \label{onevector}
\end{equation}
The right-hand side of this explicit expression makes no
mentioning of $\z_1$ anymore. It might seem at first sight that we
have shifted the problem, because it does depend on the reference
section $\z_0$. However, we have argued before that (\ref{kform})
does not depend on the choice of such a reference section, whence
our newly defined concept makes sense.

The explicit formula (\ref{onevector}) further shows that $i_{\pmbsig} \omega$ is well defined
as a $(k-1)$-form, in the sense of Definition~3. The first term on the right identifies its
associated form on $Sec(\ov \pi)$, whereas the second term, upon swapping the first two
arguments, reveals that $(i_{\pmbsig} \omega)_0=i_{\pmbsig}\omega_0$. As for local computations,
it follows from the definition (\ref{kformvector}) that $e^0(\pmbsig)=0$, whereas the
${\pmb e}^\a$ act on $\pmbsig$ simply as duals of $Sec(\ov \pi)$.

\begin{dfn} The exterior derivative of $\omega$, denoted by $d\omega$ is defined by
\begin{eqnarray}
d\omega(\z_1,\ldots,\z_{k+1}) &=& \sum_{i=1}^{k+1}(-1)^{i-1}\l(\z_i)\Big(
\omega(\z_1,\ldots,\hat{\z_i},\ldots,\z_{k+1})\Big) \nonumber \\
&& \mbox{} + \sum_{1\leq i<j\leq k+1}(-1)^{i+j}\omega([\z_i,\z_j],\z_1,\ldots,\hat{\z_i},
\ldots,\hat{\z_j},\ldots,\z_{k+1}). \label{d}
\end{eqnarray}
\end{dfn}
Note first that we are making use of definition~5 in the second term on the right, because
the bracket of two affine sections is a vector section. It is fairly obvious that $d\omega$
is skew-symmetric in all its arguments. To justify the definition, however, we should be
able to identify an operator $(d\omega)_0$ and a related $(k+1)$-form ${\bf d}\pmbom$
on $Sec({\ov \pi})$, such that $d\omega$ satisfies an expression of type (\ref{kform}).
We of course have an exterior derivative at our disposal for the $k$-form $\pmbom$ on
$Sec({\ov \pi})$ which we denote by $d$ also. We know that $d\pmbom$ has a property of type
(\ref{d}) (or can be defined that way), with vector sections replacing affine sections and
$\rho$ as anchor map instead of $\l$.

\begin{dfn} For $\omega_0: Sec(\pi)\times Sec({\ov \pi}) \times\cdots\times
Sec({\ov \pi})\rightarrow \cinfty{M}$, we define $d\omega_0$, an operator of the same
type, but depending on one more vector section, by
\begin{eqnarray}
\lefteqn{d\omega_0(\z,\pmbzeta{2},\ldots,\pmbzeta{k+1}) = \l(\z)\Big(\pmbom(\pmbzeta{2},
\ldots,\pmbzeta{k+1})\Big) } \nonumber \\
&& \mbox{} + \sum_{i=2}^{k+1}(-1)^{i-1}\r(\pmbzeta{i})\Big(\omega_0(\z,\pmbzeta{2},\ldots,
\hat{\pmbzeta{i}},\ldots,\pmbzeta{k+1})\Big) \nonumber \\
&& \mbox{} + \sum_{j=2}^{k+1}(-1)^{j+1}\pmbom([\z,\pmbzeta{j}],\pmbzeta{2},\ldots,
\hat{\pmbzeta{j}},\ldots,\pmbzeta{k+1}) \nonumber \\
&& \mbox{} - \sum_{2\leq i<j\leq k+1}(-1)^{i+j}\omega_0(\z,[\pmbzeta{i},\pmbzeta{j}],
\pmbzeta{2},\ldots,\hat{\pmbzeta{i}},\ldots,\hat{\pmbzeta{j}},\ldots,\pmbzeta{k+1}).
\label{d0}
\end{eqnarray}
\end{dfn}

This expression may look rather exotic at first, but it is obtained by formally copying
the definition (\ref{d}) and writing in that process either $\l$ or $\r$, and either
$\omega_0$ or $\pmbom$, in such a way that every term in the right-hand side has a
proper meaning. There are two important observations to be made here. First of all, the
required linearity of $d\omega_0$ in its vector arguments relies on the properties (\ref{compat1})
and (\ref{compat3}) of our Lie algebroid bracket. Secondly, replacing the affine section $\z$ in
the definition by $\z+\pmbsig$, we find:
\begin{equation}
d\omega_0(\z+\pmbsig,\pmbzeta{2},\ldots,\pmbzeta{k+1}) = d\omega_0(\z,\pmbzeta{2},\ldots,\pmbzeta{k+1})
+ d\pmbom(\pmbsig,\pmbzeta{2},\ldots,\pmbzeta{k+1}). \label{d0d}
\end{equation}
We thus know what to expect for the decomposition (\ref{kform}) of $d\omega$ and this is confirmed
by the following result.

\begin{prop} We have $(d\omega)_0=d\omega_0$ and ${\bf d}\pmbom=d\pmbom$.
\label{prop2}
\end{prop}
The proof, which involves a rather technical but straightforward
calculation, is given in the appendix.

It is of some interest to work out some simple cases in detail.
For a function $f \in \cinfty M$, $df$ is defined by
\begin{equation}
df(\z)=\l(\z)(f)=\l(\z_0)(f) + \r(\pmbzeta{})(f), \label{df}
\end{equation}
from which we learn that $(df)_0=df$ (as expected) and ${\bf df}(\pmbsig)=\r(\pmbsig)(f)$.
If $\theta$ is a 1-form, the defining relation (\ref{d}) for its exterior derivative reads
\begin{equation}
d\theta(\z_1,\z_2)=\l(\z_1)(\theta(\z_2))-\l(\z_2)(\theta(\z_1))-\pmbtheta([\z_1,\z_2]).
\label{dtheta}
\end{equation}
Introducing an arbitrary reference section  $\z_0$, it is easy to verify that this can be
rewritten as
\begin{equation}
d\theta(\z_1,\z_2)=(d\theta)_0(\z_0,\pmbzeta{2}) - (d\theta)_0(\z_0,\pmbzeta{1})
+ {\bf d}\pmbtheta(\pmbzeta{1},\pmbzeta{2}), \label{dtheta2}
\end{equation}
where ${\bf d}\pmbtheta= d\pmbtheta$ and
\begin{equation}
(d\theta)_0(\z,\pmbsig)=\l(\z)(\pmbtheta(\pmbsig))-\r(\pmbsig)(\theta(\z))
- \pmbtheta([\z,\pmbsig]). \label{dtheta0}
\end{equation}
This is in perfect agreement with the results of Proposition~3 and Definition~7.

Concerning derivation properties, it is trivial to verify that for the product of functions:
$d(fg)=fdg+gdf$. Also, from (\ref{dtheta}) applied to $f\theta$ we get:
\begin{eqnarray*}
d(f\theta)(\z_1, \z_2) &=& \l(\z_1)(f\theta(\z_2)) - \l(\z_2)(f\theta(\z_1))
-f\pmbtheta([\z_1,\z_2]) \\
&=& df(\z_1) \theta(\z_2) - df(\z_2)\theta(\z_1) + f\Big(\l(\z_1)(\theta(\z_2))
- \l(\z_2)(\theta(\z_1)) - \pmbtheta([\z_1,\z_2])\Big),
\end{eqnarray*}
from which we conclude that
\begin{equation}
d(f\theta) = f\,d\theta + df\wedge\theta. \label{dftheta}
\end{equation}
Recalling now the general statements about derivations we made before, we can conclude that
there exists a unique derivation $\hat{d}$ on $\forms{}$, of degree 1, which coincides with
our $d$ on functions and 1-forms. If we can show that $\hat{d}\omega=d\omega$ for an arbitrary
$\omega\in\forms{}$, we will know that the operator $d$ defined by (\ref{d}) is a derivation.
To this end, let us introduce
\[
\hat{d}_\z = [i_\z,\hat{d}],
\]
which, as commutator of two derivations, is itself a derivation of degree 0 on $\forms{}$. We
extend the action of $\hat{d}_\z$ to $Sec(\pi)$ `by duality'. That is to say, for $\eta\in Sec(\pi)$,
$\hat{d}_\z\eta$ is defined by requiring that for all $\theta\in\forms{1}$:
\begin{equation}
\langle \hat{d}_\z\eta,\theta\rangle = \hat{d}_\z(\theta(\eta))-\hat{d}_\z\theta(\eta). \label{dual}
\end{equation}
It is easy to see that $\hat{d}_\z\eta$ is skew-symmetric in $\z$ and $\eta$, so that it makes sense
to introduce a bracket notation for it: $[\z,\eta]^{\scriptscriptstyle\wedge}=\hat{d}_\z\eta$.
We now recall a result proved in \cite{MaCaSaII} which, although stated there in an entirely different
context, has a quite universal validity.

\begin{lemma} Given a derivation $\hat{d}$ of degree 1, and introducing $\hat{d}_\z$ and
$[\z,\eta]^{\scriptscriptstyle\wedge}$ as above, we have for all $\omega\in\forms{k}$:
\begin{eqnarray*}
\hat{d}\omega(\z_1,\ldots,\z_{k+1}) &=& \sum_{i=1}^{k+1}(-1)^{i-1}\hat{d}_{\z_i}\Big(
\omega(\z_1,\ldots,\hat{\z_i},\ldots,\z_{k+1})\Big) \\
&& \mbox{} + \sum_{1\leq i<j\leq k+1}(-1)^{i+j}\omega([\z_i,\z_j]^{\scriptscriptstyle\wedge},
\z_1,\ldots,\hat{\z_i},\ldots,\hat{\z_j},\ldots,\z_{k+1}).
\end{eqnarray*}
\end{lemma}

To show that $\hat{d}\omega=d\omega$ now, it suffices to verify that $\hat{d}_\z f=\l(\z)(f)$
on functions, and that the bracket $[\,,\,]^{\scriptscriptstyle\wedge}$ coincides with the Lie
algebroid bracket. We have $\hat{d}_\z f=i_\z\hat{d}f=i_\z df=\l(\z)(f)$, and for all $\theta\in
\forms{1}$
\begin{eqnarray*}
\langle [\z,\eta]^{\scriptscriptstyle\wedge},\theta\rangle &=& \langle\hat{d}_\z\eta,\theta\rangle
= \l(\z)(\theta(\eta)) - \hat{d}\theta(\z,\eta) - \hat{d}i_\z\theta(\eta) \\
&=& \l(\z)(\theta(\eta)) - \l(\eta)(\theta(\z)) - d\theta(\z,\eta) = \pmbtheta([\z,\eta]),
\end{eqnarray*}
from which the desired result follows.

We now reach the main question which is about the relationship between $d^2$ and the compatibility
requirements in the definition of an affine Lie algebroid. To appreciate the meaning of the
following lemma, we take a step back and assume now that the bracket $[\z_i,\z_j]$ figuring in
the definition (\ref{d}) of $d$ satisfies the `Leibniz-type property' (\ref{compat1}) with respect to
the module structure of $Sec(\pi)$ (and the
resulting property (\ref{compat3})), but no further compatibility or Lie algebra conditions a priori.
Remember that the property (\ref{compat1}) of the bracket was necessary to make sure that $d\omega$ is a
form in the first place.

\begin{lemma} For all $\omega\in\forms{k}$ and $\z_i\in Sec(\pi)$, we have
\begin{eqnarray}
\lefteqn{d^2\omega(\z_1,\ldots,\z_{k+2})= } \nonumber \\
&& \hspace*{-5mm} \sum_{1\leq i<j \leq k+2} (-1)^{i+j}\Big(\r([\z_i,\z_j]) - [\l(\z_i),\l(\z_j)]\Big)
\Big(\omega(\z_1,\ldots,\hat{\z_i},\ldots,\hat{\z_j},\ldots,\z_{k+2})\Big) \nonumber \\
&& + \sum_{1\leq i<j<l \leq k+2} (-1)^{i+j+l}\omega
({\textstyle\sum_{i,j,l}} [\z_i,[\z_j,\z_l]],\z_1,\,.\,,\hat{\z_i},
\,.\,,\hat{\z_j},\,.\,, \hat{\z_l},\,.\,,\z_{k+2}), \label{dd}
\end{eqnarray}
(where the smaller summation sign of course refers again to a cyclic sum over the three indices involved).
\end{lemma}

In fact this lemma, with suitable adaptations, also has a rather universal validity.
For completeness, we include a proof in Appendix~A.

\begin{prop} The exterior derivative has the property $d^2=0$ if and only if the bracket further
satisfies the Jacobi identity (\ref{jacobi2}) (or equivalently (\ref{formaljacobi})).
\label{prop4}
\end{prop}

\proof\ \ If (\ref{compat1}) and (\ref{jacobi2}) hold true, we also have (\ref{formaljacobi}) and
know from previous considerations that (\ref{compat2}) then holds as well.
The above lemma this way trivially implies $d^2=0$.
For the converse, we observe that $d^2f=0$, for $f \in\cinfty M$, implies (\ref{compat4}), from which
it subsequently follows that $d^2\theta=0$, with $\theta\in\forms{1}$, implies (\ref{jacobi2}). \qed

It remains now to list coordinate expressions for the basic exterior derivatives.
Let $(t,x^i)$ as before be coordinates on $M$. For their exterior derivatives we obtain
the following: for all $\z\in Sec(\pi)$,
\[
dt(\z)=\l(\z)(t)=1, \qquad dx^i(\z)=\l(\z)(x^i),
\]
from which it follows that ${\bf dt}={\bf 0}$ and ${\bf dx^i}(\pmbsig)=\r(\pmbsig)(x^i)$
(and of course $(dt)_0=dt,\ (dx^i)_0=dx^i$). In terms of the general representation (\ref{localtheta})
of a 1-form, we thus have:
\begin{eqnarray}
dt &=& e^0, \label{dt} \\
dx^i &=& \l^i e^0 + \r^i_{\a}{\pmb e}^{\a}. \label{dx}
\end{eqnarray}
Obviously, we have $de^0=0$.
We further calculate, making use, for example, of (\ref{compom0},\ref{compom}), the general formula
(\ref{dtheta}) and the coordinate expressions (\ref{strucfunc}), that
\begin{equation}
d{\pmb e}^{\a} = - C_{\b}^{\a}e^0 \wedge{\pmb e}^{\b} - \onehalf C_{\b \g}^{\a}{\pmb e}^{\b}
\wedge{\pmb e}^{\g}. \label{de}
\end{equation}
It is instructive to verify that expressing the properties $d^2{\pmb e}^\a=0$ and $d^2x^i=0$ is indeed
equivalent to the requirements (\ref{fjaccoord},\ref{jaccoord}) and
(\ref{compat2coord},\ref{liealghomcoord}), respectively.

To complete the picture of basic derivations on $\forms{}$, we have a closer look at the analogue
of the classical Lie derivative.

\begin{dfn}
For every $\z \in Sec(\pi)$, the derivation $d_\z$ of degree zero is defined as
\begin{equation}
d_\z = [i_\z,d] = i_{\z} \circ d + d \circ i_{\z}. \label{liederiv}
\end{equation}
\end{dfn}

So, since $d_\z$ is defined as a commutator of derivations, we know that it will itself be a derivation
of degree zero:
\begin{equation}
d_\z(\omega\wedge\mu) = d_\z\omega\wedge \mu +
\omega\wedge d_\z\mu. \label{liederivation}
\end{equation}
Likewise, we can rely on proofs similar to those in the standard theory to conclude that the following
commutator properties will hold true:
\begin{equation}
[d_\z,i_{\eta}] = i_{[\z,\eta]}, \qquad [d_\z,d]= 0, \qquad
[d_\z,d_\eta] = d_{[\z,\eta]}. \label{commutators}
\end{equation}
Note, however, that a Lie-type derivation with respect to a vector section turns up in the last
property, and this is indeed well defined also as: $d_{\pmbsig}=[i_{\pmbsig},d]$. It is further natural
to extend the action of $d_\z$ to $Sec(\pi)$ by duality, i.e.\ to require that a property of type
(\ref{dual}) holds true. It then follows, as expected, that for $\eta,\z\in Sec(\pi)$,
\begin{equation}
d_\z \eta = [\z,\eta]. \label{liederivdual}
\end{equation}
As a result of such an extension, $d_\z$ has Leibniz-type properties also with respect to the evaluation
of forms on the appropriate number of affine (or vector) sections; the following property, which could be
verified by a direct computation from the definition of $d_\z$, thus becomes self-evident:
\begin{equation}
d_\eta\omega (\z_1, \ldots, \z_k) = \l(\eta)\Big( \omega( \z_1
,\ldots, \z_k )\Big) + \sum_{j=1}^k (-1)^j
 \omega ([\eta,\z_j], \z_1,\ldots, \hat{\z_j},\ldots, \z_k). \label{selfduality}
\end{equation}

In the interest of doing computations, we list the Lie-type derivatives of functions $f\in\cinfty{M}$
and the local basis of 1-forms. For $\z=e_0+\z^\a{\pmb e}_\a$,
\[
d_\z f=\l(\z)(f), \qquad d_\z e^0=0, \qquad d_\z {\pmb e}^\a = C^\a_\b\z^\b e^0 - C^\a_\b{\pmb e}^\b
+ C^\a_{\b\g}\z^\g {\pmb e}^\b + d\z^\a.
\]

For future developments, it may be of some interest, finally, to list what the two composing parts
of $d_\z\omega$ are, in the sense of the defining relation (\ref{kform}) of forms.

\begin{prop} For $\omega\in\forms{k}$, we have
\begin{eqnarray}
(d_\z\omega)_0(\z_0, \pmbzeta 1, \ldots,\pmbzeta{k-1}) &=&
\l(\z)\Big( \omega_0 (\z_0, \pmbzeta 1, \ldots,\pmbzeta{k-1})\Big)
- \pmbom([\z,\z_0],\pmbzeta 1, \ldots,\pmbzeta{k-1}) \nonumber \\
&& \mbox{} + \sum_{j=1}^{k-1}(-1)^j \omega_0(\z_0,[\z,\z_j],\pmbzeta 1, \ldots,
\hat{\pmbzeta j},\ldots,\pmbzeta{k-1}), \label{liederiv0} \\
{\bf d}_{\pmbzeta{}}\pmbom(\pmbzeta 1,\ldots, \pmbzeta{k}) &=&
\l(\z)\Big( \pmbom(\pmbzeta 1, \ldots,\pmbzeta{k})\Big) \nonumber \\
&& \mbox{} + \sum_{j=1}^{k-1}(-1)^j \pmbom([\z,\z_j],\pmbzeta 1, \ldots, \hat{\pmbzeta j} \ldots,
\pmbzeta{k}). \label{liederivbold}
\end{eqnarray}
\end{prop}
These are exactly the sort of expressions one expects. The proof is a matter of a direct
computation, starting from the formula (\ref{selfduality}) and using the decompositions
(\ref{kform}) and (\ref{onevector}). It further requires manipulations of double sums of the same
nature as those in Appendix~A.

\section{$\l$-admissible curves and dynamics}

As we expressed in the introduction, the model of affine Lie
algebroids we are developing should in the first place offer an
environment in which one can accomodate the time-dependent
Lagrange-type equations (\ref{lagreq}). At present, we wish to
look in more detail at the geometric nature of the more general
dynamical systems, which we call pseudo-second-order equations,
and are those described by differential equations of the form
(\ref{dynsyst}). For this purpose in fact, we do not need the full
machinery of algebroids: it suffices to assume that $E$ is an
affine bundle over $M$ and $\l:E\rightarrow J^1M$ an affine bundle
map over the identity.

\begin{dfn}
A curve $\psi$ in $E$, which is a section of $\tau\circ\pi$, is said to be $\l$-admissible,
if $\l\circ\psi=j^1(\pi\circ\psi)$.
\end{dfn}
One could say that $\psi$ is the $\l$-prolongation of a curve in $M$. In coordinates, we have
\[
\psi:t\mapsto (t,x^i(t),y^\a(t)), \qquad \mbox{with}\quad \dot{x}^i(t)=\l^i(t,x(t))
+ \r^i_\a(t,x(t))\,y^\a(t).
\]
Note in passing that, not unexpectedly, one can characterize
$\l$-admissi\-bi\-li\-ty via a concept of contact forms: putting $\Theta^i = \l^*\theta^i$,
where the $\theta^i$ are the contact forms on $J^1M$,
we have that $\psi$ is a $\l$-admissible curve in $E$ if and only if $\psi^* \Theta^i = 0$.

Pseudo-second-order equation fields on $E$ are vector fields whose
integral curves all are $\l$-admissible curves. As in the standard
theory of \sode s on a tangent bundle or first jet bundle,
however, there is a simple direct characterization of such vector
fields.

\begin{dfn}
$\G\in\vectorfields{E}$ is a pseudo-second-order equation field if
\[
T\pi \circ \G  = i \circ \l,
\]
where $i$ is the injection of $J^1M$ into $TM$.
\end{dfn}

Clearly, in coordinates, a pseudo-\sode\ is of the form
\begin{equation}
\G=\fpd{}{t} + (\l^i(t,x) + \r_{\a}^i(t,x) y^{\a}) \fpd{}{x^i} + f^\a(t,x,y)\fpd{}{y^\a},
\label{presode}
\end{equation}
for some functions $f^\a$, and it is obvious that all its integral curves will be $\l$-admissible.

The following diagram visualizes the notions of $\l$-admissible
curves and pseudo-\sode s.

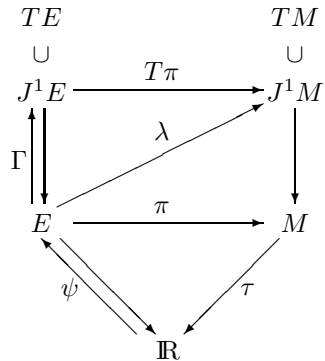
\setlength{\savelen}{\unitlength}
\setlength{\unitlength}{12pt}

\begin{picture}(12,12)(-7,-1)

\put(6.5,10){\mybox{$TE$}}
\put(14.5,10){\mybox{$TM$}}

\put(6.5,8.8){\mybox{$\cup$}}
\put(14.5,8.8){\mybox{$\cup$}}

\put(7.5,7.7){\vector(1,0){6}}  

\put(6.5,7.7){\mybox{$\JE$}}
\put(14.5,7.7){\mybox{$\JM$}}

\put(7.5,3.5){\vector(1,0){6}}  
\put(10.3,4){\mybox{$\pi$}}
\put(7,4){\vector(2,1){6.5}}   

\put(6.5,3.5){\mybox{$E$}}
\put(14.5,3.5){\mybox{$M$}}
\put(10.3,6.4){\mybox{$\l$}}

\put(6.6,7.1){\vector(0,-1){3}}  
\put(6.2,4.1){\vector(0,1){3}}  
\put(14.5,7.1){\vector(0,-1){3}}  

\put(10.3,8.3){\mybox{$T\pi$}}
\put(5.8,5.5){\mybox{$\Gamma$}}
\put(13,1.4){\mybox{$\tau$}}

\put(10.5,-0.5){\mybox{$\R$}}
\put(9.5,0){\vector(-1,1){3}} 
\put(7.1,3){\vector(1,-1){3}} 
\put(14,3){\vector(-1,-1){3}} 
\put(7.4,1.4){\mybox{$\psi$}}
\end{picture}

\setlength{\unitlength}{\savelen}

An important point now, however, is that there is a natural way of interpreting the
vector field $\G$ as section of a different bundle.

From the above definition, it is clear that a pseudo-\sode\ is
actually a section of $\otaupi: \JE \rightarrow E$, with the
additional property that for all $p\in E$, $T\pi|_{\JE} (\G(p))
=\l(p)$. An equivalent way of saying the same thing, by definition
of the concept of a pullback bundle, is that $(p,\G(p))$ is a
point of $\l^*\JE$, with $\JE$ regarded as fibred over $\JM$ via
$T\pi|_{\JE}$. From now on, we will write $\le$ for $\l^*\JE$, and
denote its two projections as indicated in the following diagram:

\begin{picture}(10,4.2)(-6,.6)
\put(4,1.4){\vector(1,0){2}}
\put(4,4){\vector(1,0){2}}
\put(3.2,3.5){\vector(0,-1){1.6}}
\put(6.8,3.5){\vector(0,-1){1.6}}
\put(3.2,1.4){\makebox(0,0){$E$}}
\put(3.2,4){\makebox(0,0){$\le$}}
\put(6.8,1.4){\makebox(0,0){$\JM$}}
\put(6.8,4){\makebox(0,0){$\JE$}}
\put(5,1){\makebox(0,0){$\l$}}
\put(5,4.35){\makebox(0,0){$\l^1$}}
\put(2.8,2.7){\makebox(0,0){$\pi_2$}}
\put(7.4,2.7){\makebox(0,0){$T\pi$}}
\end{picture}

If we finally put $\pi_1=\otaupi\circ\l^1$, there is yet another
way of expressing the characterization of a pseudo-\sode. Indeed,
from the trivial observation that
$\otaupi(\G(p))=p=\pi_2((p,\G(p)))$, it follows that a
pseudo-\sode\ $\G$ can be regarded also as a section of the bundle
$\pi_1:\le\rightarrow E$, with the property that
$\pi_2\circ\G=\pi_1\circ\G$.

The various spaces and projections, described in this discussion, are depicted in the
diagram of the next section. This diagram immediately suggests the following question:
if we put more structure into the scheme by assuming now again that $\pi:E\rightarrow M$
carries an affine Lie algebroid structure, is it possible to prolong this structure
to the bundle $\pi_1:\le\rightarrow E$, in such a way of course that $\l^1$ becomes the
anchor map of the induced affine Lie algebroid?

\section{Prolongation of affine Lie algebroids}

We shall now look in more detail at the bundle $\pi_1:\le\rightarrow E$. Its
total space is the manifold
\[
\le = \l^*\JE =\{ (q,Z) \in E \times \JE \mid \l(q) = T\pi|_{\JE}(Z) \},
\]
but the fibration we want to focus on is not one of the projections which define $\le$,
but rather the map $\pi_1=\otaupi\circ\l^1$. As such, we are looking at
an affine bundle, modelled on the vector bundle ${\ov \pi}_1:\lv\rightarrow E$, with total space
\[
\lv = \{ ({\pmb v},{\pmb V}) \in V \times VE \mid \r({\pmb v}) = T\pi|_{VE}({\pmb V}) \}.
\]
The affine bundles involved, and their underlying vector bundles are illustrated below.

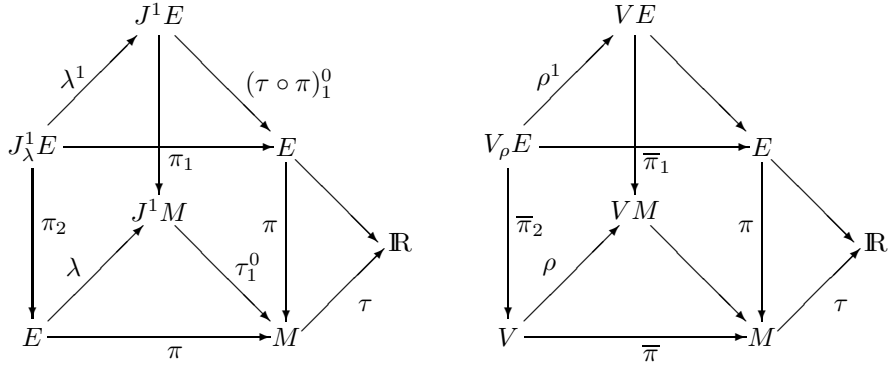
\setlength{\unitlength}{12pt}

\begin{picture}(40,13)(0,2)

\put(7.5,9.5){\vector(1,0){6.6}}  
\put(11,13){\vector(1,-1){3}}  
\put(7,10.2){\vector(1,1){2.8}}   

\put(6.5,9.5){\mybox{$\le$}}
\put(14.5,9.5){\mybox{$E$}}
\put(10.5,13.7){\mybox{$\JE$}}

\put(7,3.5){\vector(1,0){7}}  
\put(11,7){\vector(1,-1){3}}  
\put(7,4){\vector(1,1){3}}   

\put(6.5,3.5){\mybox{$E$}}
\put(14.5,3.5){\mybox{$M$}}
\put(10.5,7.5){\mybox{$\JM$}}

\put(6.5,8.8){\vector(0,-1){4.7}}  
\put(10.5,13){\vector(0,-1){5}}  
\put(14.5,9){\vector(0,-1){5}}  

\put(11,3){\mybox{$\pi$}}
\put(7.8,5.7){\mybox{$\l$}}
\put(13.3,5.7){\mybox{$\otau$}}

\put(11.2,9){\mybox{$\pi_1$}}
\put(7.8,11.7){\mybox{$\l^1$}}
\put(14.6,11.5){\mybox{$\otaupi$}}

\put(7.2,7){\mybox{$\pi_2$}}
\put(14,7){\mybox{$\pi$}}

\put(15,3.7){\vector(1,1){2.6}}  
\put(17,4.5){\mybox{$\tau$}}
\put(18.1,6.5){\mybox{$\R$}}
\put(14.8,9.1){\vector(1,-1){2.6}} 


\put(22.5,9.5){\vector(1,0){6.6}}  
\put(26,13){\vector(1,-1){3}}  
\put(22,10.2){\vector(1,1){2.8}}   

\put(21.5,9.5){\mybox{$\lv$}}
\put(29.5,9.5){\mybox{$E$}}
\put(25.5,13.7){\mybox{$VE$}}
\put(22,3.5){\vector(1,0){7}}  
\put(26,7){\vector(1,-1){3}}  
\put(22,4){\vector(1,1){3}}   

\put(21.5,3.5){\mybox{$V$}}
\put(29.5,3.5){\mybox{$M$}}
\put(25.5,7.5){\mybox{$VM$}}

\put(21.5,8.8){\vector(0,-1){4.7}}  
\put(25.5,13){\vector(0,-1){5}}  
\put(29.5,9){\vector(0,-1){5}}  

\put(26,3){\mybox{$\ov \pi$}}
\put(22.8,5.7){\mybox{$\r$}}

\put(26.2,9){\mybox{$\ov{\pi}_1$}}
\put(22.8,11.7){\mybox{$\r^1$}}

\put(22.2,7){\mybox{$\ov{\pi}_2$}}
\put(29,7){\mybox{$\pi$}}

\put(30,3.7){\vector(1,1){2.6}}  
\put(32,4.5){\mybox{$\tau$}}
\put(33.1,6.5){\mybox{$\R$}}
\put(29.8,9.1){\vector(1,-1){2.6}} 

\end{picture}
\setlength{\unitlength}{\savelen}

A section $Z$ of $\pi_1$ is completely determined once we know the maps $\pi_2\circ Z:E\rightarrow E$
and $\l^1\circ Z:E\rightarrow \JE$. Likewise, vector sections $\pmb Z$ of ${\ov \pi}_1$ are
determined by ${\ov \pi}_2\circ{\pmb Z}$ and $\r^1\circ {\pmb Z}$. For example, let $e\in E$ be a point
with coordinates $(t,x^i,y^\a)$, so that $(t,x^i)$ are the coordinates of $\pi(e)\in M$ and $e$
has the representation $e=e_0+y^\a \base{\a}$. If then $Z$ is a section of $\pi_1$, we will have:
\begin{eqnarray*}
\pi_2\circ Z: && (t,x,y)\longmapsto (t,x,z^\a(t,x,y)), \\
\l^1\circ Z: && (t,x,y)\longmapsto \left. \left(\fpd{}{t}+(\l^i+\rho^i_\a z^\a)\fpd{}{x^i}
+ Z^\a\fpd{}{y^\a}\right)\right|_{\,e},
\end{eqnarray*}
and determining $Z$ in coordinates of course amounts to assigning the functions $(z^\a,Z^\a)$ on $E$.

It is worthwhile looking at the representation of such a $Z$ with respect to suitably selected
local sections of $\pi_1$ and ${\ov \pi}_1$, which will exhibit the affine structure of $\pi_1$
and are adapted to the basis which was selected to coordinatize $E$. To this end, we introduce
two sets of local sections $\baseX \a$ and $\baseV \a$ of ${\ov \pi}_1$ which will span $Sec({\ov \pi}_1)$,
and select a zero section ${\mathcal E}_0$ as follows.
The $\baseV \a$ span `vertical sections' and are determined by:
${\ov \pi}_2\circ \baseV{\a} ={\pmb 0}$, while for $e\in E_m$ we let $\r^1\circ\baseV{\a}(e)$ be
the tangent vector to the curve $s\mapsto e+s\,\base{\a}(m)$ in $E_m$. Verticality is an intrinsic property
whereas, as usual, there is no intrinsic notion of horizontality. The determination of the $\baseX \a$
and ${\mathcal E}_0$ will therefore rely on pure coordinate arguments.
For the projection onto $V$ we put ${\ov \pi}_2\circ \baseX{\a} = \base{\a}\circ\pi$ and then,
fixing $\r^1\circ\baseX{\a}$ (as a vector field on $E$) further requires making a prescription for
the vertical components, which we simply take to be zero. Similarly, for the choice of a zero section,
we could take any vector field on $E$ which projects under $T\pi$ onto
$\l(e_0)\in\vectorfields{M}$ (and as such defines also a section of $\pi_1$), but we will fix it also by
taking the vertical components to be zero. Thus we have:
\begin{equation}
\baseX \a (e)= \left(\base{\a}(\pi(e)),\left.\r^i_{\a}(t,x)\fpd{}{x^i} \right|_e\right)
\qquad  \baseV \a(e) = \left({\pmb 0}(\pi(e)),\left.\fpd{}{y^{\a}}\right|_e\right), \label{prolongvector}
\end{equation}
and
\begin{equation}
{\mathcal E}_0(e) = \left(e_0(\pi(e)),\left.\left(\fpd{}{t} + \l^i(t,x)\fpd{}{x^i}\right)\right|_e\right).
\label{prolongzero}
\end{equation}
The general section $Z$ of $\pi_1$ then has the local representation:
\begin{equation}
Z= {\mathcal E}_0 + z^\a(t,x,y)\baseX{\a} + Z^\a(t,x,y)\baseV{\a}. \label{prolongZ}
\end{equation}

Note that pseudo-\sode s, as discussed in the previous section,
are precisely those sections $\G$ of $\pi_1$, for which
$z^\a(t,x,y)=y^\a$.

Let now $E$ be equipped with an affine Lie algebroid structure. To be in line with the notations we
used in Definition~1, we will from now on also write $\l^1(Z)$ instead of $\l^1\circ Z$, and likewise
for the $\pi_2$-projection and the corresponding projections of vector sections.

We wish to establish that there is an induced Lie algebroid structure on the affine bundle $\pi_1$.
To this end, following the scheme of Definition~1, we have to identify a bracket on ${\ov \pi}_1$ and
an action of affine sections on vector sections, such that all the necessary requirements are met.
The idea is to define such brackets by requiring roughly that its two projections are determined by the
known brackets of the projected sections. But there are some technical complications which we will
address now.

For ${\pmb Z}_1,{\pmb Z}_2\in Sec({\ov \pi}_1)$, a preliminary observation is that the Lie bracket
of their image under $\r^1$ (which gives rise to vector fields on $E$), belongs to the image of $\r^1$.
A coordinate calculation can confirm this. Putting
\[
\r^1({\pmb Z}_i)= z^\a_i\r^j_\a\fpd{}{x^j} + Z^\a_i\fpd{}{y^\a},
\]
we have
\begin{eqnarray*}
[\r^1({\pmb Z}_1),\r^1({\pmb Z}_2)] &=& \Big(\r^1({\pmb Z}_1)(z^\a_2)-\r^1({\pmb Z}_2)(z^\a_1)\Big)
\r^j_\a\fpd{}{x^j} \\
&& \mbox{} + \Big(z^\a_2\r^1({\pmb Z}_1)(\r^j_\a) - z^\a_1\r^1({\pmb Z}_2)(\r^j_\a)\Big)\fpd{}{x^j}
+ \cdots\fpd{}{y^\a}.
\end{eqnarray*}
The first term on the right manifestly belongs to the image of $\r^1$, whereas the last term is
irrelevant for that purpose. The middle term can be rewritten as
\[
z^\a_2 z^\b_1\left(\r^j_\b\fpd{\r^i_\a}{x^j} - \r^j_\a\fpd{\r^i_\b}{x^j}\right)\fpd{}{x^i},
\]
which is seen to belong to the image of $\r^1$ in view of the property (\ref{liealghomcoord}).
It is therefore natural to impose right away that the bracket $[\cdot,\cdot]^1$ under construction,
which of course is required to be skew-symmetric and $\R$-bilinear, should satisfy
\begin{equation}
\r^1\left([{\pmb Z}_1,{\pmb Z}_2]^1\right) = [\r^1({\pmb Z}_1),\r^1({\pmb Z}_2)].
\label{prolongbracket1}
\end{equation}
This will have for consequence that for $F_i\in\cinfty{E}$,
\begin{eqnarray*}
\r^1\left([F_1\,{\pmb Z}_1,F_2\,{\pmb Z}_2]^1\right) &=& \\
&& \hspace*{-3cm} F_1F_2\,[\r^1({\pmb Z}_1),\r^1({\pmb Z}_2)]
+ F_1\,\r^1({\pmb Z}_1)(F_2)\,\r^1({\pmb Z}_2)-F_2\,\r^1({\pmb Z}_2)(F_1)\,\r^1({\pmb Z}_1) .
\end{eqnarray*}
It remains then to make sure that the projection under ${\ov \pi}_2$ can be specified in a
compatible way. The above coordinate calculation to some extent illustrates how one should proceed.
If we apply $T\pi$ to the preceding equality, we get (pointwise)
\begin{eqnarray*}
\lefteqn{T\pi\left(\r^1\left([F_1\,{\pmb Z}_1,F_2\,{\pmb Z}_2]^1\right)\right) =
(F_1F_2)\,T\pi\left([\r^1({\pmb Z}_1),\r^1({\pmb Z}_2)]\right)} \\
&& \mbox{} + F_1\,\r^1({\pmb Z}_1)(F_2)\,\r({\ov \pi}_2({\pmb Z}_2))
- F_2\,\r^1({\pmb Z}_2)(F_1)\,\r({\ov \pi}_2({\pmb Z}_1)).
\end{eqnarray*}
In general, the $\r({\ov \pi}_2({\pmb Z}_i))$ are vector fields along $\pi$ for which there is no
standard Lie bracket available. If the ${\pmb Z}_i$ are projectable, however, meaning that
there exist $\pmbzeta{i}\in Sec({\ov \pi})$ such that
${\ov \pi}_2\circ {\pmb Z}_i=\pmbzeta{i}\circ\pi$, the vector fields $\r^1({\pmb Z}_i)$ on $E$
are $\pi$-related to the vector fields
$\r(\pmbzeta{i})$ on $M$. Hence, the corresponding brackets are also $\pi$-related, meaning that
for projectable ${\pmb Z}_i$, we can put
\begin{equation}
{\ov \pi}_2\left([{\pmb Z}_1,{\pmb Z}_2]^1\right) = [{\ov \pi}_2({\pmb Z}_1),{\ov \pi}_2({\pmb Z}_2)],
\label{prolongbracket2}
\end{equation}
and then the property (\ref{liealghom}) (which in coordinates gives (\ref{liealghomcoord})) ensures that
\[
T\pi\left(\r^1\left([{\pmb Z}_1,{\pmb Z}_2]^1\right)\right) =
\r\circ{\ov \pi}_2\left([{\pmb Z}_1,{\pmb Z}_2]^1\right)
\]
as it should. The expression for $T\pi\left(\r^1\left([F_1\,{\pmb Z}_1,F_2\,{\pmb Z}_2]^1\right)\right)$
further shows that the ${\ov \pi}_2$ and $\r^1$ projections of the bracket under construction
will still match up if for projectable
${\pmb Z}_i$ and for any $F_i\in\cinfty{E}$, we define
\begin{eqnarray}
\lefteqn{{\ov \pi}_2\left([F_1\,{\pmb Z}_1,F_2\,{\pmb Z}_2]^1\right) =
F_1F_2\,[{\ov \pi}_2({\pmb Z}_1),{\ov \pi}_2({\pmb Z}_2)]} \nonumber \\
&& \mbox{} + F_1\,\r^1({\pmb Z}_1)(F_2)\,{\ov \pi}_2({\pmb Z}_2)
- F_2\,\r^1({\pmb Z}_2)(F_1)\,{\ov \pi}_2({\pmb Z}_1). \label{prolcompat1}
\end{eqnarray}
It then follows that
\begin{equation}
[F_1\,{\pmb Z}_1,F_2\,{\pmb Z}_2]^1=F_1F_2\,[{\pmb Z}_1,{\pmb Z}_2]^1
+ F_1\,\r^1({\pmb Z}_1)(F_2)\,{\pmb Z}_2 - F_2\,\r^1({\pmb Z}_2)(F_1)\,{\pmb Z}_1, \label{prolcompat2}
\end{equation}
since both sides have the same ${\ov \pi}_2$ and $\r^1$ projections.

The final point to observe now is that sections of ${\ov \pi}_1$ (locally) are finitely generated, over the
ring $\cinfty{E}$, by projectable sections. Hence, the defining relations (\ref{prolongbracket1})
and (\ref{prolcompat1}) are sufficient to define the bracket $[\cdot,\cdot]^1$ on vector sections.
The property (\ref{prolcompat2}) will hold by extension for all vector sections and the bracket
will satisfy the Jacobi identity as a result of the Jacobi identity of the Lie algebroid bracket
we start from and the same identity for vector fields on $E$.

To define the action of $Z\in Sec(\pi_1)$ on ${\pmb V}\in Sec({\ov \pi}_1)$, we proceed in exactly
the same manner. First, one easily verifies that the Lie bracket of $\l^1(Z)$ and $\r^1({\pmb V})$
belongs to the image of $\r^1$, this time in view of the properties (\ref{compat2coord}) and
(\ref{liealghomcoord}). Hence, it makes sense to put
\begin{equation}
\r^1\left([Z,{\pmb V}]^1\right) = [\l^1(Z),\r^1({\pmb V})], \label{prolongbracket3}
\end{equation}
and we of course require the bracket $[\cdot,\cdot]^1$ to have linearity properties of the kind
of (\ref{linearity}).
For projectable sections, we can put
\begin{equation}
{\ov \pi}_2\left([Z,{\pmb V}]^1\right) = [\pi_2(Z),{\ov \pi}_2({\pmb V})],
\label{prolongbracket4}
\end{equation}
and be assured of consistency with the projection (\ref{prolongbracket3}). Next, still for
projectable $Z$ and ${\pmb V}$ and for any $F\in\cinfty{E}$, we define
\begin{equation}
{\ov \pi}_2\left([Z,F\,{\pmb V}]^1\right) = F\,[\pi_2(Z),{\ov \pi}_2({\pmb V})]
+ \l^1(Z)(F)\,{\ov \pi}_2({\pmb V}). \label{prolcompat3}
\end{equation}
Sections of $\pi_1$ can be written as a projectable zero section, plus a linear combination of
projectable vector sections with coefficients in $\cinfty{E}$. In combination with the earlier
arguments for vector sections, we are again led to the conclusion that the requirements
(\ref{prolongbracket3}) and (\ref{prolcompat3}) are sufficient to define $[Z,{\pmb V}]^1$ for
arbitrary $Z\in Sec(\pi_1)$ and ${\pmb V}\in Sec({\ov \pi}_1)$ and that we will have the
property
\begin{equation}
[Z,F\,{\pmb V}]^1 = F\,[Z,{\pmb V}]^1 + \l^1(Z)(F)\,{\pmb V}. \label{prolcompat4}
\end{equation}
The final requirement of type (\ref{formaljacobi}) then also easily follows, which concludes
the construction of the prolonged affine Lie algebroid.

For computational purposes, it remains to list the brackets of the local sections which are used
in the general representation of a section of $\pi_1$ as in (\ref{prolongZ}). We have
\begin{eqnarray*}
&& [ {\mathcal E}_0, \baseX \a ]^1 = C_\a^\b \baseX \b, \qquad
[ {\mathcal E}_0, \baseV \a ]^1 = {\pmb 0}, \\
&& [\baseX \a, \baseX \b]^1 = C_{\a\b}^\g \baseX \g, \quad
[\baseX \a, \baseV \b]^1 = {\pmb 0},\quad [\baseV \a, \baseV \b]^1 = {\pmb 0}.
\end{eqnarray*}
It is perhaps worthwhile to repeat hereby that the two projections have to be looked at
to verify these statements, although of course they are bound to match up
if our new bracket has been defined consistently. Thus we have, for example:
\begin{eqnarray*}
{\ov \pi}_2\left([{\mathcal E}_0, \baseX \a]^1\right) &=& [\pi_2({\mathcal E}_0),
{\ov \pi}_2(\baseX \a)] = [e_0,\base{\a}] = C_\a^\b \base{\b} , \\
\r^1\left([{\mathcal E}_0, \baseX \a]^1\right) &=& [\l^1({\mathcal E}_0),
\r^1(\baseX \a)] = \left[\fpd{}{t}+\l^i\fpd{}{x^i},\r^j_\a\fpd{}{x^j}\right] = C^\b_\a\,
\r^1(\baseX \b),
\end{eqnarray*}
where (\ref{compat2coord}) has been used again in the last line.

\section{Discussion and outlook for future work}

The form of equations (\ref{lagreq}), which we claim to be the appropriate generalization of
Lagrangian systems on Lie algebroids to a situation where explicit time-dependence is involved,
has brought us to the introduction of the new concept of Lie algebroids on affine bundles which
are fibred over $\R$. More precisely, the first guidance for developping this concept was
provided by the conditions (\ref{oldceq}) and (\ref{newceq}) which the various functions
appearing in (\ref{lagreq}) have to satisfy. Ultimately, of course, we want to arrive at an
intrinsic geometrical construction of such Lagrangian systems. There are many aspects to be
explored yet, but we have sufficiently paved the way already to be able to predict what the
outcome of subsequent studies will bring.

One of us has shown \cite{Mart} that the prolongation of a Lie
algebroid (in the standard situation of vector bundles) provides a
platform where there exist analogues of the intrinsic structures
living on a tangent bundle and these in turn give rise to an
intrinsic definition of Lagrangian systems via Poincar\'e-Cartan
type forms. This is the reason why we were keen to verify
immediately that the same notion of prolongation exists in our
affine set-up. There is little doubt now that we will find
intrinsic objects on such a prolonged affine Lie algebroid, which
are analogues of what is known to give rise to an intrinsic
definition of time-dependent Lagrangian systems on the first jet
bundle of a manifold fibred over $\R$. But there is more to it.
Even when there is no Lagrangian for the dynamics under
consideration and we are, in other words, talking about
pseudo-\sode s on a Lie algebroid, we expect to be able to develop
a machinery of associated non-linear and linear connections, which
again is analogous to the standard theory of connections
associated to \sode s on a tangent bundle or first jet bundle. In
fact, a paper on these issues for the case of Lie algebroids on a
vector bundle is in preparation.

One of the features we examined in this paper as a kind of test for the relevance and internal
consistency of the generalized notion of Lie algebroids, was the existence of an associated
coboundary operator $d$. But of course, there are still other interesting properties which
standard Lie algebroids are known to exhibit. Let us briefly highlight another one here
and show that it also survives our generalization, namely the existence of an associated
Poisson structure. Specifically, we want to establish that there exists a canonically
defined Poisson structure on the extended dual $E^\dag$.

Sections of $\pi$ (respectively $\ov{\pi}$) can be identified with
linear functions on $E^\dag$ (respectively $V^*$). Explicitly, if
$\z\in Sec(\pi)$, we consider the function $\hat{\z}\in
\cinfty{E^\dag}$ defined by: for each $p\in E^\dag$, $p\in
E^\dag_m$ say, $\hat{\z}(p)=p(\z_m)$. Likewise, if $\pmbsig\in
Sec(\ov{\pi})$, we denote by $\hat{\pmbsig}\in\cinfty{V^*}$ the
function defined by $\hat{\pmbsig}({\boldsymbol p})={\boldsymbol
p}(\pmbsig_m)$, where ${\boldsymbol p}\in V^*_m$. In coordinates,
if $\z=e_0+\z^\a \base{\a}$ and $p\in E^\dag$ has coordinates
$(t,x^i,p_0,p_\a)$, then $\hat{\z}(p)=p_0+p_\a\z^\a(t,x)$ and
similarly for $\hat{\pmbsig}$. Now, for any two sections $\z,\e\in
Sec(\pi)$, we define the function $\{\hat{\z},\hat{\e}\}$ on
$E^\dag$ by
\begin{equation}
\{\hat{\z},\hat{\e}\}(p) = \widehat{[\z,\e]}({\boldsymbol p}),
\label{poisson1}
\end{equation}
whereby we recall that $[\z,\e]$ is a section of $\ov{\pi}$ and
for $p\in E^\dag_m$, ${\boldsymbol p}$ is the associated element
of $V^*_m$. If further $f,g$ are functions on $M$ and we make no
notational distinction for their pullback to $E^\dag$, internal consistency of (\ref{poisson1})
for the action of $Sec(\ov{\pi})$ on $Sec(\pi)$ requires that we further put
\begin{equation}
\{\hat{\z},f\}=-\{f,\hat{\z}\}= \l(\z)(f), \qquad \{f,g\}=0.
\label{poisson2}
\end{equation}
The construction then uniquely extends to a skew-symmetric, $\R$-bilinear
bracket operation on $E^\dag$ with the required derivation property. This bracket satisfies the
Jacobi identity as a result of the Lie algebroid Jacobi identity (\ref{jacobi2}).

The brackets for the coordinate functions on $E^\dag$ are found to be:
\begin{align*}
& \{t,t\}=0 &\quad &\{t,x^i\}=0 &\quad &\{x^i,x^j\}=0 \\
& \{p_0,t\}=1 &\quad & \{p_0,x^i\}=\l^i  &\quad &\{p_0,p_\beta\}=C^\gamma_{\beta}p_\gamma \\
& \{p_\alpha,t\}=0 &\quad &\{p_\alpha,x^i\}=\rho^i_\alpha &\quad
&\{p_\alpha,p_\beta\}=C^\gamma_{\alpha\beta}p_\gamma.
\end{align*}
There is an interesting observation to be made here. Recall that $E^\dag$ is actually a vector
bundle and note now that the bracket we have constructed preserves the subset of functions
on $E^\dag$ which are linear in the fibre coordinates. As a result, we know that there is an
induced Lie algebroid structure on the bundle $(E^\dag)^*\rightarrow M$. There are many new
insights to be gained from approaching the subject of a Lie algebroid structure on the affine
bundle $E\rightarrow M$ from this angle; that is to say, by regarding $E\rightarrow M$ as an
affine subbundle of $(E^\dag)^*\rightarrow M$ and taking an appropriate Lie algebroid structure
on $(E^\dag)^*\rightarrow M$ as the starting point. Also this will be the subject of a forthcoming
paper.


\section*{Appendix:\ \ Technical proofs}

The start for proving Proposition~3 is the defining relation (\ref{d}) of the exterior derivative,
in which we make use of the decomposition (\ref{kform}) in the first term and (\ref{onevector})
in the second. We first obtain,
\begin{eqnarray*}
\lefteqn{d\omega(\z_1,\ldots,\z_{k+1}) } \\
&& = \sum_{i=1}^{k+1}(-1)^{i-1}(\l(\z_0)+\r(\pmbzeta{i}))
\Big(\sum_{j=1}^{i-1}(-1)^{j-1}\omega_0(\z_0,\pmbzeta 1,\ldots,\hat{\pmbzeta j}, \ldots,
\hat{\pmbzeta i},\ldots,\pmbzeta{k+1}) \\
&& \mbox{} + \sum_{j=i+1}^{k+1}(-1)^{j}\omega_0(\z_0,\pmbzeta 1,\ldots,\hat{\pmbzeta i}, \ldots,
\hat{\pmbzeta j},\ldots,\pmbzeta{k+1}) +
\pmbom(\pmbzeta 1,\ldots,\hat{\pmbzeta i},\ldots,\pmbzeta{k+1})\Big) \\
&& \mbox{} + \sum_{1\leq i <j\leq k+1} (-1)^{i+j}\Big( \pmbom([\z_i,\z_j],\pmbzeta 1,\ldots,
\hat{\pmbzeta i},\ldots,\hat{\pmbzeta j},\ldots,\pmbzeta{k+1}) \\
&& \mbox{} + \sum_{l=1}^{i-1}(-1)^{l}\omega_0(\z_0,[\z_i,\z_j],\pmbzeta 1,\ldots,\hat{\pmbzeta l},
\ldots,\hat{\pmbzeta i},\ldots,\hat{\pmbzeta j},\ldots,\pmbzeta{k+1}) \\
&& \mbox{} + \sum_{l=i+1}^{j-1}(-1)^{l-1}\omega_0(\z_0,[\z_i,\z_j],\pmbzeta 1,\ldots,\hat{\pmbzeta i},
\ldots,\hat{\pmbzeta l},\ldots,\hat{\pmbzeta j},\ldots,\pmbzeta{k+1}) \\
&& \mbox{} + \sum_{l=j+1}^{k+1}(-1)^{l}\omega_0(\z_0,[\z_i,\z_j],\pmbzeta 1,\ldots,\hat{\pmbzeta i},
\ldots,\hat{\pmbzeta j},\ldots,\hat{\pmbzeta l},\ldots,\pmbzeta{k+1})\Big),
\end{eqnarray*}
and now perform a number of manipulations on multiple sums. Interchanging the order of summation in the
first line, we have $\sum_{i=1}^{k+1}\sum_{j=1}^{i-1}=\sum_{j=1}^{k}\sum_{i=j+1}^{k+1}$. Interchanging
subsequently the names of the indices $i$ and $j$, the term involving $\l(\z_0)$ of the first line cancels
the similar one in the second line. The last three lines involve triple sums, which can be rearranged
as follows. The first triple sum, with suitable interchanges of the order of summation, becomes:
\[
\sum_{i=1}^{k}\sum_{j=i+1}^{k+1}\sum_{l=1}^{i-1}= \sum_{i=1}^{k}\sum_{l=1}^{i-1}\sum_{j=i+1}^{k+1}
=\sum_{l=1}^{k-1}\sum_{i=l+1}^{k}\sum_{j=i+1}^{k+1}=\sum_{1\leq l<i<j\leq k+1}.
\] For the second one, we have
\[
\sum_{i=1}^{k}\sum_{j=i+1}^{k+1}\sum_{l=i+1}^{j-1}=\sum_{i=1}^{k}\sum_{l=i+1}^{k}\sum_{j=l+1}^{k+1}
=\sum_{1\leq i<l<j\leq k+1}.
\]
The last one can directly be written as $\sum_{1\leq i<j<l\leq k+1}$. Changing names of indices to make
all triple sums look alike, we thus far arrive at the result:
\begin{eqnarray*}
\lefteqn{d\omega(\z_1,\ldots,\z_{k+1}) } \\
&& = \sum_{1\leq i<j\leq k+1}(-1)^{i+j}\r(\pmbzeta j -\pmbzeta i)\Big(\omega_0(\z_0,\pmbzeta 1,
\ldots,\hat{\pmbzeta i}, \ldots,\hat{\pmbzeta j},\ldots,\pmbzeta{k+1})\Big) \\
&& \mbox{} + \sum_{i=1}^{k+1}(-1)^{i-1}(\l(\z_0)+\r(\pmbzeta{i}))
\Big(\pmbom(\pmbzeta 1,\ldots,\hat{\pmbzeta i},\ldots,\pmbzeta{k+1})\Big) \\
&& \mbox{} + \sum_{1\leq i <j\leq k+1} (-1)^{i+j}\, \pmbom([\z_0,\pmbzeta{j}]-[\z_0,\pmbzeta{i}]
+[\pmbzeta{i},\pmbzeta{j}],\pmbzeta 1,\ldots,
\hat{\pmbzeta i},\ldots,\hat{\pmbzeta j},\ldots,\pmbzeta{k+1}) \\
&& \mbox{} + \sum_{1\leq i<j<l\leq k+1}(-1)^{i+j+l}
\omega_0(\z_0,[\z_i,\z_j]+[\z_j,\z_l]+[\z_l,\z_i],\pmbzeta 1,\,.\,,\hat{\pmbzeta i},
\,.\,,\hat{\pmbzeta j},\,.\,,\hat{\pmbzeta l},\,.\,,\pmbzeta{k+1}).
\end{eqnarray*}
It is clear now that the terms which do not involve $\z_0$ combine exactly to
$d\pmbom(\pmbzeta{1},\ldots,\pmbzeta{k+1})$. What remains is
\begin{eqnarray*}
&& \sum_{i=1}^{k+1}(-1)^{i-1}\l(\z_0)
\Big(\pmbom(\pmbzeta 1,\ldots,\hat{\pmbzeta i},\ldots,\pmbzeta{k+1})\Big) \\
&& \mbox{} + \sum_{1\leq i<j\leq k+1}(-1)^{i+j}\r(\pmbzeta j -\pmbzeta i)\Big(\omega_0(\z_0,\pmbzeta 1,
\ldots,\hat{\pmbzeta i}, \ldots,\hat{\pmbzeta j},\ldots,\pmbzeta{k+1})\Big) \\
&& \mbox{} + \sum_{1\leq i <j\leq k+1} (-1)^{i+j}\, \pmbom([\z_0,\pmbzeta{j}]-[\z_0,\pmbzeta{i}],\pmbzeta 1,
\ldots,\hat{\pmbzeta i},\ldots,\hat{\pmbzeta j},\ldots,\pmbzeta{k+1}) \\
&& \mbox{} + \sum_{1\leq i<j<l\leq k+1}(-1)^{i+j+l}
\omega_0(\z_0,[\z_i,\z_j]+[\z_j,\z_l]+[\z_l,\z_i],\pmbzeta 1,\,.\,,\hat{\pmbzeta i},
\,.\,,\hat{\pmbzeta j},\,.\,,\hat{\pmbzeta l},\,.\,,\pmbzeta{k+1}),
\end{eqnarray*}
and should be compared to $\sum_{i=1}^{k+1}(-1)^{i-1}d\omega_0(\z_0,\pmbzeta 1,\ldots,\hat{\pmbzeta{i}},
\ldots,\pmbzeta{k+1})$, with $d\omega_0$ as defined in (\ref{d0}). It is obvious that the first three
lines in the computation of $(d\omega)^0$ are exactly the ones we have in the above expression. The last
term in (\ref{d0}) gives rise to triple sums of the form
$\sum_{l=1}^{k+1}(-1)^{l-1}\sum_{i<j<l}(-1)^{i+j-1}\omega_0(\z_0,[\pmbzeta{i},\pmbzeta{j}],\pmbzeta{1},\,.\,,
\pmbzeta{i},\,.\,,\pmbzeta{j},\,.\,,\pmbzeta{l},\,.\,,\pmbzeta{k+1})$ (there is a similar term with
$\sum_{i<l<j}$ and one with $\sum_{l<i<j}$). With suitable interchanges of summations, similar to what was
explicitly explained before, these three terms combine to:
\[
\sum_{1\leq i<j<l\leq k+1}(-1)^{i+j+l}
\omega_0(\z_0,[\pmbzeta{i},\pmbzeta{j}]+[\pmbzeta{j},\pmbzeta{l}]+[\pmbzeta{l},\pmbzeta{i}],
\pmbzeta 1,\,.\,,\hat{\pmbzeta i},\,.\,,\hat{\pmbzeta j},\,.\,,\hat{\pmbzeta l},\,.\,,\pmbzeta{k+1}).
\]
The proof now becomes complete if we observe that:
\[
[\z_i,\z_j]+[\z_j,\z_l]+[\z_l,\z_i]=
[\pmbzeta{i},\pmbzeta{j}]+[\pmbzeta{j},\pmbzeta{l}]+[\pmbzeta{l},\pmbzeta{i}].
\]

We next turn to the proof of Lemma~2.

If $\omega$ is a $k$-form, then $d^2\omega$ is a $(k+2)$-form with
\begin{eqnarray}
d^2\omega(\z_1,\ldots,\z_{k+2}) &=&
\sum_{i=1}^{k+2}(-1)^{i-1}\l(\z_i)\Big(
d\omega(\z_1,\ldots,\hat{\z_i},\ldots,\z_{k+2})\Big) \nonumber \\
&& \mbox{} + \sum_{1\leq i<j\leq k+2}(-1)^{i+j}d\omega([\z_i,\z_j],\z_1,\ldots,\hat{\z_i},
\ldots,\hat{\z_j},\ldots,\z_{k+2}). \label{d2}
\end{eqnarray}
If we plug in the definition of $d\omega$, the first term on the right
will further decompose into two parts, one involving double
and the other involving triple sums. Based on our experience with such combinatorics in
the preceding proof, we can right away conclude that the first line of the right-hand side
of (\ref{d2}) equals:
\begin{eqnarray}
&& \hspace*{-8mm}\sum_{1\leq i<j\leq k+2}(-1)^{i+j}\Big(\l(\z_j)\l(\z_i)-\l(\z_i)\l(\z_j)\Big)
\Big(\omega(\z_1,\ldots,\hat{\z_i},\ldots,\hat{\z_j},\ldots,\z_{k+2})\Big) \nonumber \\
&& \hspace*{-8mm}\mbox{} + \sum_{1\leq i<j<l\leq k+2} (-1)^{i+j+l-1}{\textstyle\sum_{i,j,l}}\left\{
\l(\z_i)\Big(\omega([\z_j,\z_l]\right\},\z_1,\,.\,,\hat{\z_i},\,.\,,\hat{\z_j},\,.\,,
\hat{\z_l},\,.\,,\z_{k+2})\Big) , \label{d2first}
\end{eqnarray}
where the smaller summation sign, as before, refers to a cyclic sum, the range of which is delimited
by the curly brackets.
For the second term on the right in (\ref{d2}), we have to remember that the first argument
is a vector section. Using the defining relation (\ref{kformvector}), applied to $d\omega$,
we obtain:
\begin{eqnarray}
\lefteqn{d\omega(\pmbsig,\z_1,\ldots,\z_{k}) = \r(\pmbsig)\Big( \omega
(\z_1,\ldots,\z_{k}) \Big)} \nonumber \\
&& \mbox{} + \sum_{i=1}^{k}(-1)^{i}\l(\z_i)\Big(\omega(\pmbsig, \z_1,\ldots,
\hat{\z_i},\ldots,\z_{k})\Big) \nonumber \\
&& \mbox{} + \sum_{j=1}^{k}(-1)^{j}\omega([\pmbsig,\z_j],\z_1,\ldots,\hat{\z_j},
\ldots,\z_{k}) \nonumber \\
&& \mbox{} + \sum_{1\leq i<j\leq k}(-1)^{i+j}\omega([\z_i,\z_j],\pmbsig,\z_1,\ldots,\hat{\z_i},
\ldots,\hat{\z_j},\ldots,\z_{k}). \label{onevectord}
\end{eqnarray}
The last line here has two vector arguments, but this is consistent with the application of
definition~5 to a form of type $i_{\pmbzeta{}}\omega$. We look at the effect of each of these four
terms, when inserted in the second sum of (\ref{d2}).
The first one simply gives:
\begin{equation}
\sum_{1\leq<i<j\leq k+2} (-1)^{i+j}\r([\z_i,\z_j])\Big( \omega
(\z_1,\ldots,\hat{\z_i},\ldots, \hat{\z_j}, \ldots, \z_{k+2})\Big). \label{rhopart}
\end{equation}
The second one is easily seen to give rise to terms which cancel exactly the second sum in
(\ref{d2first}). The third term of (\ref{onevectord}) gives rise to expressions involving
double brackets, which combine to:
\begin{equation}
\sum_{1\leq i<j<l \leq k+2} (-1)^{i+j+l} \omega({\textstyle \sum_{i,j,l}}[[\z_i,\z_j],\z_l],
\z_1,\,.\,,\hat{\z_i},\,.\,,\hat{\z_j},\,.\,,\hat{\z_l},\,.\,,\z_{k+2}). \label{jacobipart}
\end{equation}
The fourth term of (\ref{onevectord}) finally creates terms which involve two double sums,
and in each of the summands the first two arguments of $\omega$
are brackets. One has to look at all possible orderings, six in total, of the four different
indices involved, but when the same procedure is applied to shuffle the order of summations
suitably around and rename indices where appropriate, one easily finds that the six
terms cancel each other two by two in view of the skew-symmetry of $\omega$. What we are left
with in the end is the first term of (\ref{d2first}), (\ref{rhopart}) and (\ref{jacobipart}):
they precisely combine to the statement in Lemma~2.

{\small 

}  

\end{document}